\documentclass{article}              % DO NOT DELETE THIS LINE

\usepackage{amsmath}
\usepackage{amssymb}
\usepackage{cite}
\newif\ifpdf
\ifx\pdfoutput\undefined
   \pdffalse
\else
   \pdfoutput=1
   \pdftrue
\fi
 
\ifpdf
 \usepackage[pdftex]{graphicx}
  \RequirePackage{hyperref}
  \PassOptionsToPackage{pdftex,bookmarksopen,bookmarksnumbered}{hyperref}
\else
  \usepackage{epsfig}
\fi

\newcommand{\eqr}[1]{\eqref{#1}}
\newcommand{\figr}[1]{Figure~\ref{#1}}

\newcommand{\eps}{{\epsilon}}
\newcommand{\bet}{{\beta}}

\newcommand{\ubn}{u_{\beta_{\scriptscriptstyle n}}}

\newcommand{\fbn}{f_{\beta_{\scriptscriptstyle n}}}
\newcommand{\betn}{\itr{\beta}{n}}

\newcommand{\betnpo}{\itr{\beta}{n+1}}

\newcommand{\norm}[1]{\left\|#1\right\|}

\newcommand{\itr}[2]{{#1_{{#2}}}}
\newcommand{\seq}[1]{{\itr{#1}{0},\itr{#1}{1},\itr{#1}{2},\dots }}
\newcommand{\Fix}{\ensuremath{\operatorname{Fix}\,}}

\newcommand{\spl}{{S_+}}
\newcommand{\FT}[1]{\mathcal{F}{#1}}
\newcommand{\IFT}[1]{\mathcal{F}^{-1}{#1}}
\newcommand{\set}[1]{{\field #1}}

\newcommand{\map}[3]{#1:\,#2\rightarrow #3\,}
\newcommand{\field}[1]{\mathbb{#1}}

\newcommand{\NN}{\mathbb{N}}

\newcommand{\ZZ}{\mathbb{Z}}

\newcommand{\rp}{\set{R}_+}

\newcommand{\vh}{\widehat{v}}
\newcommand{\T}{{\ensuremath{\mathcal T}}}

\newcommand{\asrncv}{{\T_{*}}}
\newcommand{\asr}{{T_{*}}}
\newcommand{\rasr}{V(\asr,\beta)}

\newcommand{\rasrn}{V(\asr,\betn)}
\newcommand{\rasrnpo}{V(\asr,\betnpo)}
\newcommand{\rasrncvx}{V(\asrncv,\bet)}
\newcommand{\rasrncvbn}{V(\asrncv,\betn)}
\newcommand{\rasre}{V({\widetilde{\asrncv}},\beta)}
\newcommand{\rasrne}{V({\widetilde{\asrncv}},\betn)}
\newcommand{\RB}{{R_{\scriptscriptstyle \!B}}}
\newcommand{\RA}{{R_{\scriptscriptstyle \!A}}}
\newcommand{\PB}{{P_{\scriptscriptstyle \!B}}}
\newcommand{\PA}{{P_{\scriptscriptstyle \!A}}}
\newcommand{\PC}{{P_{\scriptscriptstyle \!C}}}
\newcommand{\RC}{{R_{\scriptscriptstyle \!C}}}
\newcommand{\RM}{{R_{\scriptscriptstyle \!M}}}
\newcommand{\PM}{{P_{\scriptscriptstyle \!M}}}
\newcommand{\PS}{{P_{\scriptscriptstyle \!S}}}

\newcommand{\RSP}{{R_{\scriptscriptstyle \!{S_+}}}}
\newcommand{\RS}{{R_{\scriptscriptstyle \!{S}}}}
\newcommand{\PSP}{{P_{\scriptscriptstyle \!{S_+}}}}
\newcommand{\esp}{E_{\scriptscriptstyle \!\spl}}
\newcommand{\espn}{E_{\scriptscriptstyle \!\spl}{(\itr{x}{n})}}
\newcommand{\psp}{P_{\scriptscriptstyle \!\spl}}

\newcommand{\ind}[1]{{{\mathcal X}_{#1}}}
\newcommand{\indCD}{{\ind{D^c}}}
\newcommand{\indD}{{\ind{D}}}

\newcommand{\LL}{\ensuremath{\mathcal L}}

\newcommand{\emp}{\ensuremath{\mathrm{\mbox{\rm\O}}}}

\newcommand{\supp}{\ensuremath{\mbox{supp}}}
\newcommand{\closu}[1]{\ensuremath{\mbox{cl}{#1}}}

\def\proof{\noindent{\it Proof}. \ignorespaces}
\def\endproof{\ensuremath{\quad \blacksquare}}

\begin{document}                  % DO NOT DELETE THIS LINE

\title{\sffamily Relaxed Averaged Alternating
  Reflections\\ for Diffraction Imaging}

\author{D.\ Russell Luke\thanks{
The Pacific Institute for the Mathematical Sciences,
Simon Fraser University, Burnaby, British Columbia V5A 1S6, Canada.
E-mail: \texttt{rluke@cecm.sfu.ca}.}
}

\date{May 11, 2004 --- Version 2.2}

\maketitle

\begin{abstract}
We report on progress in algorithms for iterative
phase retrieval.  The theory of convex optimization
is used to develop and to gain insight into counterparts for the
nonconvex problem of phase retrieval.  We propose a relaxation of
averaged alternating reflectors and determine the fixed point set of
the related operator in the convex case.  A numerical study supports 
our theoretical observations and demonstrates the
effectiveness of the algorithm compared to the current state of
the art.
\end{abstract}

\section{Introduction}
\label{sec:intro}

The phase retrieval problem is a classical inverse problem in optics 
that has received renewed interest in applications to
nonperiodic scatterers and macromolecules.
While scattering from some structures allows one explicitly
to compute the phase from magnitude measurements \cite{Hauptman},
more general classes of scatterers require the use of less direct methods.  
So called {\em iterative transform methods} pioneered by Gerchberg and Saxton
\cite{Gerch72}, and Fienup \cite{Fien82} are well established generic
techniques for iteratively recovering the phase in a variety of
settings.   Recent developments in imaging
\cite{AllenOxleyPaganin,Chapman,Cloetens1,Cloetens2,Paganin,Spence1,Spence2}, 
have placed a premium on improving
the efficiency and stability of phase retrieval algorithms.

In this work we derive a stable and fast new strategy for phase retrieval,
what we call Relaxed Averaged Alternating Reflection (RAAR),
that falls under the category of iterative transform methods \cite{Luke02a}. 
The motivation for the RAAR algorithm comes from 
recent work in which another new algorithm, the Hybrid Projection
Reflection algorithm (HPR), was presented \cite{BCL2}.
The HPR algorithm was originally concieved as a single
parameter relaxation of the well known Douglas-Rachford algorithm
applied to phase retrieval.  The HPR algorithm can also be viewed a special 
case of the three parameter {\em difference map} recently proposed by 
Elser \cite{Elser1}.  

There are two fundamental and distinct issues that accompany these algorithms.
The first is the incorporation of {\em a priori} information into the 
constraint structure of the algorithms.  The second is the choice of 
algorithm parameter values.  Regarding the first issue, it is difficult to 
overestimate the effect of the constraints on the mathematical properties
and performance of the algorithms.  There have been several studies on the choice of
constraints in applications to crystallography \cite{Millane1,Millane2,Elser2}.  
We use a simple example to illustrate how seemingly minor changes 
in the physical domain constraints can lead to algorithms that 
appear very different. This has caused some confusion in the literature which
we hope to clarify through an examination of the abstract algorithmic structures
behind the leading techniques.  The choice of parameters also
has a dramatic impact on the mathematical properties of the algorithms 
and hence performance.  Physical insight often provides the best (and only) 
basis for chosing values for the algorithm parameters, 
but this is not always available or reliable.  
In the case of the HPR algorithm, our numerical experiments have 
not provided an empirical basis upon which to make recommendations.
A more mathematically rigorous approach also appears to be difficult and
has been found in only a few very special cases.  For instance,  
in the convex setting the 
convergence properties of the HPR algorithm are known for the unrelaxed
case \cite{BCL3}.  For the relaxed HPR algorithm and the more general 
difference map a complete and mathematically rigorous analysis has 
yet to be found.  To circumvent the analytical
barriers facing the difference map and the HPR algorithm, we introduce the
RAAR algorithm, which is conceptually simple, analytically tractable
and easy to implement;  moreoever, it outperforms the current state of the art.  
While the RAAR algorithm coincides with the HPR
algorithm in a limiting case, it does not fall in the class of 
algorithms covered by Elser's difference map framework.  

A precise statement of the leading algorithms is given in  Section \ref{s:RAAR}  
In this same section we provide a terse outline of the mathematical
justification for the RAAR algorithm.  In Section \ref{results} we demonstrate
the effectiveness of the algorithm and make practical recommendations for
implementation.

\section{Phase Retrieval and Iterative Transform Algorithms}\label{s:RAAR}
\subsection{Phase retrieval}
We are interested in recovering
the scattering amplitude ${u}_*$ of a medium that has
been illuminated by an electromagnetic wave from measurements of its
spatial coherence function and other {\em a priori} information.
For the sake of concreteness, we assume that ${u}_*$ is a real-valued,
nonnegative function supported on some prescribed bounded set $D$, that is
$\LL\ni\map{{u_*}}{\ZZ^N}{\rp}$ with
$\supp(u_*)\subset D\subset \ZZ^N$.
Here $\LL$ is a Hilbert space of square integrable functions, $\ZZ^N$ 
is the domain -- in this case the {\em physical domain} --
corresponding to discrete (i.e. sampled) waves, $\rp$ is the
positive reals and $\supp(u_*)$ is the support
of $u_*$.  Writing this in terms of constraints, we have $u_*\in
S_+\subset\LL$, where $S_+$ is the set of nonnegative functions in 
$\LL$ with support on $D$.  If we require only that the functions be 
supported on $D$, we denote the corresponding constraint set by $S$.  The
sets $S$ and $\spl$ are refered to as the {\em physical domain constraints}.
The other constraint we consider comes from the the data, $m$,
which we presume consists of noisy magnitude measurements in the far field, 
thus $m$ is proportional to the modulus
of the Fourier transform of ${u}_*$.  We therefore refer to the domain 
of the image data $m$ as the {\em Fourier domain}.  In terms of 
constraint sets, we write that $u_*\in M$ where
$
M=\big\{v\in\LL\;\mid\;\;|\FT{v}| = m\big\}
$
and $\FT{v}$ denotes the discrete Fourier transform of $v$.  
We shall refer to the set $M$ as the {\em Fourier}, or {\em image domain 
constraint}.  Note that $S_+$ is a {\em convex} set, while $M$ is {\em
  nonconvex}.  It is the nonconvexity of the magnitude constraint that
does not allow us to transfer classical convergence results for the
most common algorithms to the case of phase retrieval.  For further
discussion see \cite{BCL1}.

\subsection{Iterative Transform Algorithms}
We formulate the problem of phase retrieval as a feasibility problem:
\[
~\text{find}\;\;{u}\in S_+\cap M.
\]
Iterative transform techniques are built upon combining 
projections onto the sets $S_+$ and $M$ in some fashion.
While they are seldom written as fixed-point algorithms, iterative
transform algorithms can usually be put into the form
$\itr{u}{n+1}=\T\itr{u}{n}$ where $\T$
is a generic operator in which the projections and averaging
operations are embedded (see \cite{BCL2,BCL1}).  For added control and
flexibility, one often includes a {\em relaxation} strategy parameterized
by $\beta$.  We write the relaxed operator with generic, single parameter 
relaxation strategy ${\mathcal V}$ (there can be infinitely many such strategies)
as ${\mathcal V}(\T,\bet)$.  In order effectively to
exploit relaxations for improved algorithm performance it is necessary to
understand the mathematical properties of the operator ${\mathcal V}(\T,\bet)$
-- first and foremost of these is
the characterization of the set of fixed points, $\Fix {\mathcal V}(\T,\bet)$.
We return to this issue at the end of this section.

The operators we study are built upon projectors and reflectors.  
Denote by $\PC $ an arbitrary but fixed selection, or {\em projector}, 
from the possibly multi-valued
{\em projection} onto a subset $C$ of $\LL$. Closely related is the  
corresponding {\em reflector} with respect to $C$
\[
\RC = 2\PC -I,
\]
where $I$ is the identity operator.  By definition, for every
$u\in\LL$,  $\PC ({u})$ is the midpoint between ${u}$ and $\RC ({u})$.
Specializing to our application, the projector, $\PM {u}$, of a signal
${u}\in\LL$ onto the Fourier magnitude constraint set $M$
is given by
\begin{equation}\label{PM}
\PM ({u}) = \IFT({\widehat{v}_0}),\qquad\mbox{where}\qquad
\widehat{v}_0(\xi) = \begin{cases}
m(\xi)\displaystyle{\frac{\FT{{u}}(\xi)}{|\FT{{u}}(\xi)|}}, & \text{if
$\FT{{u}}(\xi) \neq 0$;} \\
m(\xi), & \text{otherwise}~.
\end{cases}
\end{equation}
Here,  $\IFT{}$ is the discrete inverse Fourier transform and $\vh_0$
a selection from the multi-valued Fourier domain projection.
For further discussion of this projector see 
Luke {\em et al} \cite[Corollary 4.3]{Luke02a} and \cite{Luke03c}.  
We return to the issue of
multivaluedness of the magnitude projection in Section \ref{results}.
The projection of a signal ${u}\in\LL$ onto $S_+$ is single-valued
(since $\spl$ is convex), and is given by
\begin{equation}\label{PSP}
(\forall x \in \ZZ^N)\qquad
\big(\PSP({u})\big)(x) = \begin{cases}
\max\{0,{u}(x)\}, & \text{if $x\in D$;} \\
0, & \text{otherwise.}
\end{cases}
\end{equation}

One of the best known iterative transform algorithms is Fienup's
Hybrid Input-Output algorithm (HIO) \cite{Fien82}.  We use this as our
benchmark for performance.   In the present setting, HIO is given as
\begin{equation}
(\forall x\in\ZZ^N)\;\;
\itr{{u}}{n+1}(x) = \begin{cases}
\big(\PM (\itr{{u}}{n})\big)(x),&
\text{if $x\in D$ and $\big(\PM (\itr{{u}}{n})\big)(x)\geq 0$;} \\
\itr{{u}}{n}(x)-\betn\big(\PM (\itr{{u}}{n})\big)(x), &\text{otherwise.}
\end{cases}
\label{e:hio}
\end{equation}

There have been several attempts to identify the HIO algorithm with a broader
class of relaxation strategies that can be written as fixed point iterations, 
that is, in the form
$\itr{u}{n+1}={\mathcal V}(\T,\betn)\itr{u}{n}.$.
Bauschke, Combettes and Luke \cite{BCL1} proved that, when only a support 
constraint as opposed to support and nonegativity is 
applied in the physical domain,
then the HIO algorithm with $\beta=1$ corresponds to the classical 
Douglas-Rachford algorithm for which convergence results in the convex setting
are well known.  In a subsequent article Bauschke Combettes and Luke \cite{BCL2} proved that, for
physical domain support constraints only, 
the HIO algorithm corresponds to a particular relaxation of the 
Douglas-Rachford algorithm, that is 
\begin{equation}
(\forall x\in\ZZ^N)\;\;
\itr{{u}}{n+1}(x) = \begin{cases}
\big(\PM (\itr{{u}}{n})\big)(x),&
\text{if $x\in D$}\\
\itr{{u}}{n}(x)-\betn\big(\PM (\itr{{u}}{n})\big)(x), &\text{otherwise,}
\end{cases}
\label{e:hioS}
\end{equation}
is equivalent to 
\begin{equation}\label{e:hprS} 
\itr{u}{n+1} = 
\tfrac12\big(\RS(\RM+(\betn-1)\PM )+I+(1-\betn)\PM \big)(\itr{{u}}{n}).
\end{equation} 
Independent of these results, Elser \cite{Elser1} showed the correspondence 
between the HIO algorithm with only support constraints in the physical domain
and the difference map,
\begin{equation}\label{e:dms} 
\itr{u}{n+1} = \left(I + \beta \left(\PS\left((1-\gamma_2)\PM-\gamma_2 I\right)
+ \PM\left((1-\gamma_1)\PS - \gamma_1I\right)\right)\right)(\itr{u}{n}),
\end{equation}
for the case where $\gamma_1=-1$ and $\gamma_2=1/\beta$.  The correspondence 
between the difference map and the HIO algorithm does not carry over to 
the case of support and nonnegativity constraints.  The correct formulation 
of the corresponding algorithm was given in \cite[Proposition 2]{BCL2}, where
it is shown that
\begin{equation}
\itr{{u}}{n+1}=
\tfrac12\big(\RSP(\RM+(\betn-1)\PM )+I+(1-\betn)\PM \big)(\itr{{u}}{n}).
\label{e:hypocampe}
\end{equation}
is equivalent to 
\begin{equation}
(\forall x\in\ZZ^N)\;\;
\itr{{u}}{n+1}(x) = \begin{cases}
\big(\PM (\itr{{u}}{n})\big)(x),&
\text{if $x\in D$ and} \\
&\text{$\big(\RM(\itr{{u}}{n})\big)(x)
\geq (1-\betn)\big(\PM (\itr{{u}}{n})\big)(x)$;} \\[2mm]
\itr{{u}}{n}(x)-\betn\big(\PM (\itr{{u}}{n})\big)(x), &\text{otherwise.}
\end{cases}
\label{e:hpr}
\end{equation}
In \cite{BCL2} the fixed point iteration \eqr{e:hypocampe} is called the 
Hybrid Projection Reflection (HPR) algorithm, which is 
equivalent to the difference map (with $\gamma_1=-1$ and $\gamma_2=1/\beta$)
applied to support and nonnegativity constraints: 
\begin{equation}\label{e:dm} 
\itr{u}{n+1} = \left(I + \beta \left(\psp\left((1-\gamma_2)\PM-\gamma_2 I\right)
+ \PM\left((1-\gamma_1)\psp - \gamma_1I\right)\right)\right)(\itr{u}{n}).
\end{equation}

It is important to note that, while the form of 
prescriptions of projection algorithms in terms of fixed
point iterations $\itr{u}{n+1}={\mathcal V}(\T,\betn)\itr{u}{n}$ does not 
depend on the underlying constraints, this is not the case for
prescriptions of the form \eqr{e:hio}, \eqr{e:hioS}
and \eqr{e:hpr}.  As we have seen,
slight changes in the constraint sets can result in dramatic changes in the
form of algorithms when written in this way.
When written as fixed point iterations, the
effect of changing the constraint structure is seen in the mathematical 
properties of the operator rather than the form of the algorithm.

Preliminary numerical results indicate that the HPR algorithm is a
promising alternative to HIO -- HPR is more stable and, at least with simulated
data, produces higher quality images.  Detailed convergence results have been
obtained in \cite{BCL3} for the unrelaxed HPR algorithm ($\beta=1$) in a
convex setting.  At this time, however, 
there are no mathematically rigorous results proving convergence
or suggesting how to choose the relaxation parameter $\betn$ to 
improve performance.  Another drawback to the HPR algorithm is that, while
it consistently delivers higher quality solutions than HIO, it can take longer
to achieve this. The algorithm we propose next addresses both the analytical
drawbacks as well as the performance issues regarding the HPR algorithm and
the more general difference map.

The new algorithm we propose is given by the following:
given any $\itr{{u}}{0}\in\LL~$, generate the sequence $\seq{{u}}$ by
\begin{equation}\label{e:berlin}
\itr{{u}}{n+1} = \rasrncvbn \itr{{u}}{n}
\end{equation}
where
\begin{equation}\label{rasrncvx}
\rasrncvx = \beta \asrncv + (1-\bet)\PM\quad\mbox{and}\quad\asrncv=\tfrac12(\RSP\RM+I).
\end{equation}
To underscore the connection of this algorithm with the Averaged Alternating
Reflection (AAR) algorithm studied in \cite{BCL3}, we refer to 
(\ref{e:berlin}) as the
relaxed averaged alternating reflection (RAAR) algorithm.
For $\bet=1$ the RAAR, HPR, AAR, and the difference map ($\gamma_1=-1$ 
and $\gamma_2=1/\beta$) algorithms are equivalent.  For $\beta\neq 1$ 
the RAAR algorithm is fundamentally different than HPR; moreover, it
cannot be derived as a special case of the difference map \eqr{e:dm}.
The recursion \eqr{e:berlin} can be written analogously to \eqr{e:hio}
and \eqr{e:hpr}.  To see this, we proceed as in Proposition 2 of
\cite{BCL2}.  Given an arbitrary signal $v\in\LL$, let
$v^+=\max\{v,0\}$ and $v^-=\min\{v,0\}$ be its positive and negative
parts, respectively.  Then \eqr{e:berlin} can be rewritten as
\begin{equation}
\itr{{u}}{n+1} =\left(-\indCD\cdot\betn(2\PM-I) -
  \left[\indD\cdot\betn(2\PM-I)\right]^- +
  \PM\right)(\itr{{u}}{n}).\label{Berlin}
\end{equation}
There are $3$ cases to consider:  (i) If $x\in D$ and $(\RM
\itr{{u}}{n})(x)\geq 0$, then \eqr{Berlin} yields
$
\itr{{u}}{n+1} = \PM;~
$
(ii)
 if $x\in D$ and $(\RM \itr{{u}}{n})(x)<0$, then \eqr{Berlin} becomes
$
\itr{{u}}{n+1}(x)= 
\big(\left((1-2\betn)\PM+\betn I\right)(\itr{{u}}{n})\big)(x);~
$
(iii) if $t\notin D$, then \eqr{Berlin} can also be written as
$
\itr{{u}}{n+1}(x)= 
\big(\left((1-2\betn)\PM+\betn I\right)(\itr{{u}}{n})\big)(x).
$
Altogether this yields the following algorithm
\begin{equation}
(\forall x\in\ZZ^N)\qquad
\itr{{u}}{n+1}(x) = \begin{cases}
\big(\PM(\itr{{u}}{n})\big)(x),&
\text{if $x\in D$ and $\big(\RM(\itr{{u}}{n})\big)(x) \geq 0$;} \\[2mm]
\betn \itr{{u}}{n}(x)-(1-2\betn)\big(\PM(\itr{{u}}{n})\big)(x), 
&\text{otherwise}.\hfill
\end{cases}
\label{e:RAAR}
\end{equation}
We summarize the above discussion in the following proposition.
% \begin{propn}\label{relax}
\medskip

\noindent
{\sc Proposition 2.1.}  {\em Algorithm \eqref{e:RAAR} is equivalent to the recursion 
\eqr{e:berlin}.}
% \end{propn}
\medskip

The update rule in algorithm (\ref{e:RAAR}) depends on the pointwise 
sign of the reflector $\big(\RM(\itr{{u}}{n})\big)(x)$ whereas the 
update rule for Fienup's HIO
algorithm (\ref{e:hio}) depends on  the pointwise sign of the projector
$\big(\PM(\itr{{u}}{n})\big)(x)$.  The difference between the RAAR
update rule and that for  HPR (\ref{e:hpr}) is much starker.  Also
note that the ``otherwise'' action is
simply a relaxation of the conditional action in the HIO algorithm;  this is,
again, very different than the HPR algorithm.

\subsection{The RAAR algorithm: convex analysis}
To gain some insight into the behavior of the algorithm above, we
study the behavior of the convex analog to $\rasrncvx$.
Let $A$ and $B$ be two closed convex subsets of $\LL$.
Replace $\spl$ and $M$ by $A$ and $B$ respectively.  Let $E\subset A$
denote the set of points in $A$ nearest to $B$, and let $F\subset B$  
denote the set of points in $B$ nearest to $A$.  The {\em gap vector}
between $A$ and $B$,
denoted by $g\in \LL$, is defined by ${g}=P_{\closu{(B-A)}}(0)$.  Loosely
interpreted, this is a vector pointing from $E$ to $F$ with  $\|{g}\|$ 
measuring the smallest distance between  $A$ and $B$.  
For instance, if $A\cap B\neq \emp$
then ${g}=0$.  For a more precise treatment see \cite{BauBorSVA,BauBorJAT}.
The convex counterpart to \eqr{rasrncvx}, the
central operator in the RAAR algorithm, is defined by
\begin{equation}\label{rasr}
\rasr = \bet \asr+(1-\bet)\PB,\quad 0<\bet<1
\quad\mbox{where}\quad
\asr = \tfrac12(\RA \RB + I).
\end{equation}

When discussing convergence of projection-type algorithms, one must take
care to distinguish between {\em consistent} and {\em inconsistent} feasibility
problems.  In the current convex setting, consistent problems satisfy
$A\cap B\neq\emp$;  when $A\cap B=\emp$ the problem is said to be inconsistent.
Inconsistent problems are common in applications where the {\em a priori}
information represented by the constraint sets is highly idealized, particularly
in the presence of noise.  Bauschke, Combettes and Luke \cite{BCL3} show that 
the properties of the
AAR algorithm (that is, RAAR with $\bet=1$) for consistent problems are
very different from inconsistent problems.  The reason for this is that the
operator $\asr$ {\em does not have} a fixed point if $A\cap B=\emp$. 
For $0<\beta<1$ the convex instance of the RAAR algorithm avoids these 
complications by transferring questions of consistency of
the constraints to the existence of {\em nearest} points.
In other words, the RAAR operator enjoys the advantage that
$\Fix \rasr$ is independent of whether or not the
associated feasibility problem is consistent.  This is the content of the
following theorem.
	\medskip

% \begin{thm}\label{t:fix}
\noindent {\sc Theorem 2.2.}  {\em
Let $0<\beta<1$.  Then
\begin{equation}
\Fix \rasr = F - \frac{\bet}{1-\bet}{g}
\label{e:fix}
\end{equation}
where $g$ is the gap vector between $A$ and $B$ and $F\subset B$ is the set of points
in $B$ nearest to $A$.  Moreover, for every $u\in \Fix\rasr$, we have the following:
\begin{equation}\label{Miles}
(i)~ u = \PB u - \frac{\beta}{1-\beta}g;
\quad (ii)~\PB u-\PA \RB u = g;
\quad
(iii)~\PB u \in F~\mbox{ and }~\PA\PB u \in E.
\end{equation}
}
% \end{thm}
\medskip
\noindent\proof
To prove the result we must show (a) that
$F-\bet g/(1-\bet)\subset \Fix\rasr$
and (b), conversely, that $\Fix\rasr\subset F-\bet g/(1-\bet).$
The first statement (a) is proved analogously to the proof of
equation (18) of \cite{BCL3}.  In the interest of brevity, we leave
this as an exercise.

We show that $\Fix(\bet\asr+(1-\bet)\PB)\subset F
-\frac{\bet}{1-\bet}g$.  To see
this, pick any $u\in \Fix(\bet\asr+(1-\bet)\PB)$.  Let $f = \PB u$ and
$y=u-f$.  For any $b\in B$, since $B$ is a nonempty closed convex set
and $f=\PB u,$ we have
$
\langle b-\PB u,u-f \rangle \leq 0.
$
which yields
\begin{equation}\label{one}
\langle b-f,~y\rangle = \langle b-f,~u-f\rangle \leq 0.
\end{equation}
Recall that $\PA(2f-u) = \PA(2\PB u - u) = \PA\RB u$.
Together with the identity \cite[Proposition 3.3(i)]{BCL3}
\begin{equation}\label{nice id}
(\forall u \in \LL) \quad u-\asr u = \PB u - \PA\RB u
\end{equation}
equation \eqr{one} yields
\begin{equation}\label{Iron}
\PA(2f-u) = f + \asr u - u.
\end{equation}
Now $\bet\asr u  + (1-\bet)\PB u = u$ yields
\begin{equation}\label{Wine}
\asr u - u =
\frac{1-\bet}{\bet}(u-\PB u).
\end{equation}
Then \eqr{Iron} and \eqr{Wine} give
$
\PA(2f-u) = f + \frac{1-\bet}{\bet}(u-f) = f +\frac{1-\bet}{\bet}y.
~$
As above, for any $a\in A,~$ since $A$ is nonempty, closed and convex,
we have
$
\langle a-\PA(2f-u),(2f-u)-\PA(2f-u) \rangle \leq 0,
$
and hence
\begin{eqnarray}
0&\geq&\left\langle a - \left(f+\frac{1-\bet}{\bet}y\right), ~(2f - u) -
\left(f +\frac{1-\bet}{\bet}y\right) \right\rangle\nonumber\\
&=&\left\langle a - \left(f+\frac{1-\bet}{\bet}y\right),~ - y -
  \frac{1-\bet}{\bet}y \right\rangle\nonumber\\
&=&\frac{1}{\bet}\left\langle -a + f,~ y \right\rangle +
\frac{(1-\bet)}{(\bet)^2}\|y\|^2.
\label{two}
\end{eqnarray}
Now \eqr{one} and \eqr{two} yield
$
\left\langle b-a,~ y \right\rangle \leq
-\frac{1-\bet}{\bet}\|y\|^2\leq 0.~
$
Now take a sequence $\seq{a}$ in $A$ and a sequence $\seq{b}$ in
$B$ such that $\itr{g}{n}=\itr{b}{n}-\itr{a}{n}\to g$.  Then
\begin{equation}\label{teil1}
(\forall n\in\NN)\quad\left\langle \itr{g}{n},~ y \right\rangle \leq
-\frac{1-\bet}{\bet}\|y\|^2\leq 0.
\end{equation}
Taking the limit and using the Cauchy-Schwarz inequality yields
\begin{equation}\label{teil2}
\|y\|\leq\frac{\bet}{1-\bet}\|g\|.
\end{equation}
Conversely,
$
u-(\bet\asr u+(1-\bet)\PB u) =
\bet\left(f - \PA(2f - u)\right) + (1-\bet)y = 0~
$
gives
\begin{equation}\label{teil3}
\|y\| =
\frac{\bet}{1-\bet}\Big\| f-\PA(2f - u)\Big
\|\geq \frac{\bet}{1-\bet}\|g\|.
\end{equation}
Hence $\|y\| = \frac{\bet}{1-\bet}\|v\|$ and,
taking the limit in \eqr{teil1}, $y = -\frac{\bet}{1-\bet}g$,
which confirms (i).
It follows immediately that $f-\PA\RB u = g$ which proves
(ii) and, by definition,
implies that
$\PB u= f\in F$ and $\PA\PB u\in E$.  This yields (iii)
and proves \eqr{e:fix}.
\endproof

\noindent In words, regardless of whether or not $A\cap B$ is empty, as
long as there are points in $B$ that are nearest to $A$, then the RAAR
operator $\rasr$ has a set of fixed points, and these are precisely the
points in B nearest to A, translated by the scaled gap vector.
This is the starting point for the convex heuristics behind the RAAR
algorithm.  Statements about convergence and more detailed behavior of
the algorithm are beyond the scope of this work.

We conclude the mathematical analysis with some observations that motivate
the relaxation strategy we implement in Section \ref{results}.  We wish to
use the parameter $\bet$ to control the step size between successive iterates
and, as much as possible, to steer the iterates.  
Far away from the solution, it is easy to see the damping effect 
of the parameter $0<\beta<1$, which derives  
from the form of the relaxation \eqr{rasrncvx} as simply a
convex combination of the operator $\asrncv$ and the projector
onto the {\em data} $\PM$ -- the smaller the relaxation parameter $\bet$,
the closer to the data we require the iterates to stay.
It was noted in \cite{BCL2} that, regardless of the relaxation,
the HPR algorithm \eqr{e:hpr} takes significantly longer than the
HIO algorithm \eqr{e:hio} to reach a suitable neighborhood of the
solution, although, once near a solution, HPR delivers consistently
better images with greater stability and reliability than HIO.  
We show in the next section that the dampening effect of the
relaxation in the RAAR algorithm is just what is needed to control the 
initial behavior of the HPR algorithm.  

For the behavior of the algorithm near the solution, we rely on the convex 
analysis.  By \eqr{e:fix}, the relaxation
parameter $\bet$ effects the fixed points of the operator through the gap vector.
If the feasibility problem is consistent, that is, $A\cap B\neq\emp$,
then the gap vector $g=0$.  In this case, is it not clear what
effect, if any, $\beta$ will have on convergence.  On the other hand,
if the problem is inconsistent, that is, $A\cap B=\emp$,
and $g\neq 0$, then, by \eqr{e:fix}, the set of nearest points $F$
can be translated arbitrarily far away in the direction $g$ by letting
$\bet$ approach $1$ from below.
We use this to gain some control on the step size between successive 
iterates and the directions of the steps.
\medskip

% \begin{propn}\label{t:step size}
\noindent {\sc Proposition 2.3.}  {\em 
Let $\itr{u}{n}\in\LL$ satisfy $\norm{\itr{u}{n}-\ubn}<\delta$ 
where $\ubn\in\Fix\rasrn$ and $0<\betn<1$. Define 
$\itr{u}{n+1}=\rasrnpo \itr{u}{n}$ for any $0<\betnpo<1$.  
Then
\begin{equation}\label{the horror}
\norm{\itr{u}{n+1} - \left(\fbn-\frac{\betnpo}{1-\betn}g\right)} < \delta,
\quad\text{where}\quad\fbn=\PB \ubn\in F.
\end{equation}
}
% \end{propn}
\proof
% \noindent {\sc Proof.}
For any ${u}\in\LL$, by \eqr{nice id}, we have 
$\rasrnpo {u}-\rasrn {u}=
(\betnpo - \betn)\left(\PA - I\right)\RB{u},~
$
which, together with \eqr{Miles}(i), yields
\begin{equation}\label{hunger}
\ubn - \rasrnpo \ubn = \frac{\betnpo-\betn}{1-\betn}g,\quad \mbox{ or }\quad
\rasrnpo\ubn = \fbn - \frac{\betnpo}{1-\betn}g.
\end{equation}
Since $\rasrnpo$ is nonexpansive,
the result follows from \eqr{hunger}.
\endproof

\noindent While the HPR algorithm gives quite stable solutions eventually,
the above theory suggests that this stability can be improved in 
a controlled fashion.  Consider the fixed point iteration as a 
descent algorithm minimizing
some error metric (in fact, minimizing the gap distance) 
where $-g$ is the direction of descent. By \eqr{the horror} 
and the first equation in \eqr{hunger}, 
\[
\itr{u}{n+1}\approx\rasrnpo\ubn  = \ubn - \frac{\betnpo-\betn}{1-\betn}g,
\] 
thus one can use $\betnpo$ to affect steps in the direction $-g$ ranging, 
in the limit, from length $-\betn/(1-\betn)$ to
$1$ as $\betnpo$ varies from $0$ to $1$ respectively.  
The difference $\itr{u}{n+1}-\itr{u}{n}$ for the unrelaxed
algorithm ($\beta=1$) was shown in \cite{BCL3} to converge to the
negative gap vector $-g$ in the inconsistent case.
The effect of the relaxation is primarily to dampen the iteration
in the neighborhood of a solution in the case of inconsistent problems.  
To see the advantage of this, consider the nonconvex case and 
suppose that the problem is inconsistent 
(that is, the gap vector $g\neq 0$).  The only case of the HPR algorithm
for which we can say anything is the case $\beta=1$, which is the same 
as the unrelaxed RAAR (or AAR) algorithm, so we restrict the discussion
to the RAAR and AAR algorithms.  The convex analysis of the AAR algorithm
shows that, even though the gap is attained, the iterates $\itr{u}{n}$ 
continue to move in the direction $-g$ without end.  In the nonconvex 
setting, even if the true gap is attained, 
the continued progress of the iterates in the direction $-g$ could 
push the iterates away from the domain of attraction of the local solution
and into a different domain of attraction.  Thus the projections of the 
iterates, or the {\em shadows} might never converge.  This ``wandering" of
the iterates near an apparent local solution has been observed both with 
the HIO and HPR algorithms, though it is much less severe and destabilizing
with HPR than it is with HIO.  The relaxations in the RAAR algorithm can 
be used to either dampen the iterates near a local solution to slow 
drifting out of a domain of attraction, or to halt the wandering of 
the iterates altogether by holding the relaxation parameter at a 
fixed value less than $1$.

\section{Numerical Implementation}\label{results}

Our goal with the RAAR algorithm is to use dynamic relaxations to
shorten the initial ``warm-up" phase of the HPR algorithm and to stabilize 
the algorithm near a local solution. 
The algorithm we consider is
\begin{equation}\label{practical}
\itr{{u}}{n+1} \approx \rasrncvbn\itr{{u}}{n}.
\end{equation}

Before outlining our specific implementation, some remarks are in
order about the calculation of $\T$ given by \eqr{rasrncvx}.
As discussed in  \cite[Section 5.2]{Luke02a} the projection onto the magnitude
constraint $\PM$ is a numerically unstable operation due to the multivaluedness of the
projection operator.
We therefore recommend the following approximation to $\PM$ (see \cite[Eq.74]{Luke02a}):
\begin{equation}
\PM {u}\approx\nabla J_\eps
{u}= I - \IFT{}\left(\left(\frac{|\FT {u}|^2}{\left(|\FT {u}|^2+\epsilon^2\right)^{1/2}} -
m\right)\frac{|\FT {u}|^2+2\epsilon^2}{\left(|\FT {u}|^2+\epsilon^2\right)^{3/2}}\FT
{u}\right)
\label{gradE}
\end{equation}
for $~0 < \epsilon \ll 1~$, where
\begin{equation}\label{sdeps}
J_\eps({u}) =~
\frac{1}{2}\left(\|{u}\|^2 - \left\|\IFT{\vh-m}~\right\|^2\right),\quad\mbox{where}\quad
\vh~ = ~\frac{\left|\FT{{u}}\right|^2}{\left(\left|\FT{{u}}\right|^2+~\epsilon^2\right)^{1/2}}.
\end{equation}
Define
\begin{equation}\label{rasre}
\rasre = \frac{1}{2}\left(\RSP(2\nabla J_\eps - I) + I\right).
\end{equation}
Under reasonable assumptions, by the continuity of $\RSP$ and
\cite[Corollary 5.3]{Luke02a} it can be shown that
$~\nabla J_\eps({u}) \to \PM({u})~$ and $~\rasre {u} \to \rasrncvx {u}~$ as $\eps\to 0.$

\begin{figure}
\caption{\label{true}
{\small Original images and corresponding data used for
the comparison of the HIO and HPR algorithms.
The center of (a) is a $38\times 38$
pixel section of the standard Lena image, zero-padded to a
$128\times 128$ matrix. Frame (b) is the noiseless Fourier magnitude data $m$
corresponding to image (a). The same object domain support
constraint (and initial guess) of size $64\times 64$ pixels, shown in (c), is
used for each trial.}}
\begin{center}
	\ifpdf
	 (a) \includegraphics[height=3.0cm,width=3.0cm,angle=0]%
	 {pictures/rasr/lena64c_pad.jpg}
	 (b) \includegraphics[height=3.0cm,width=3.0cm,angle=0]%
	 {pictures/rasr/lena64c_m.jpg}
	 (c) \includegraphics[height=3.0cm,width=3.0cm,angle=0]%
	 {pictures/rasr/lena64c_S.jpg}
	\else
	 (a) \includegraphics[height=3.0cm,width=3.0cm,angle=0]%
	 {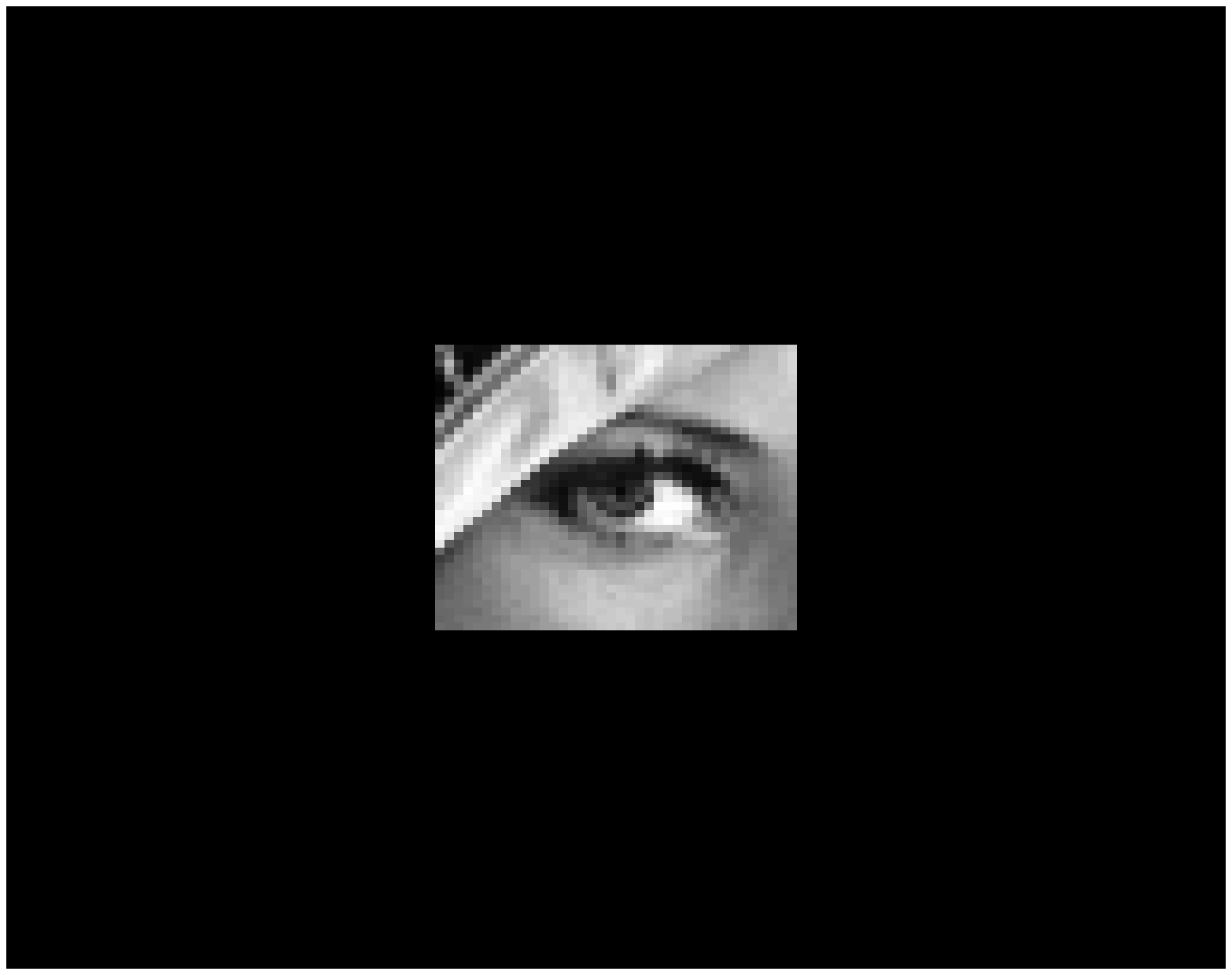}
	 (b) \includegraphics[height=3.0cm,width=3.0cm,angle=0]%
	 {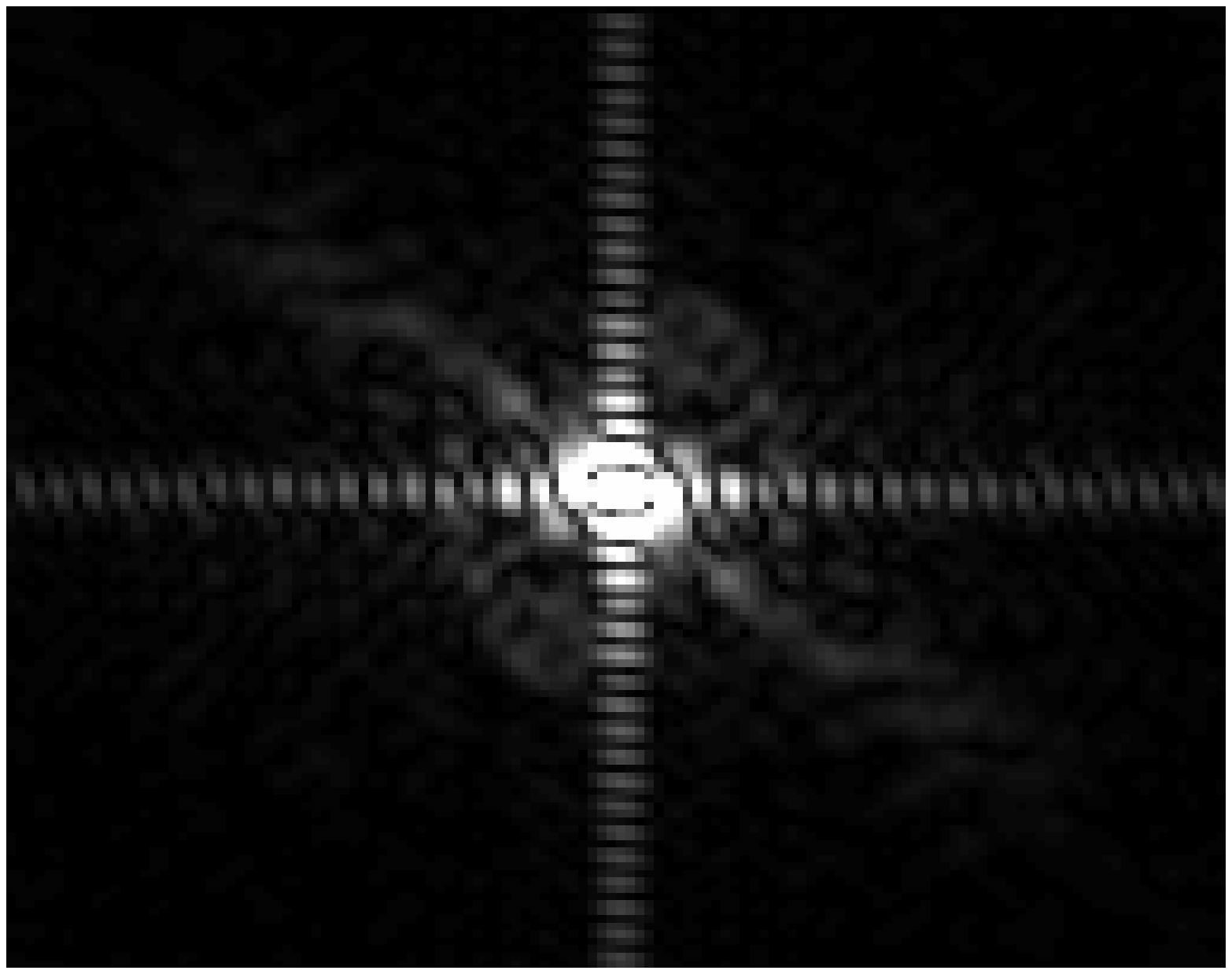}
	 (c) \includegraphics[height=3.0cm,width=3.0cm,angle=0]%
	 {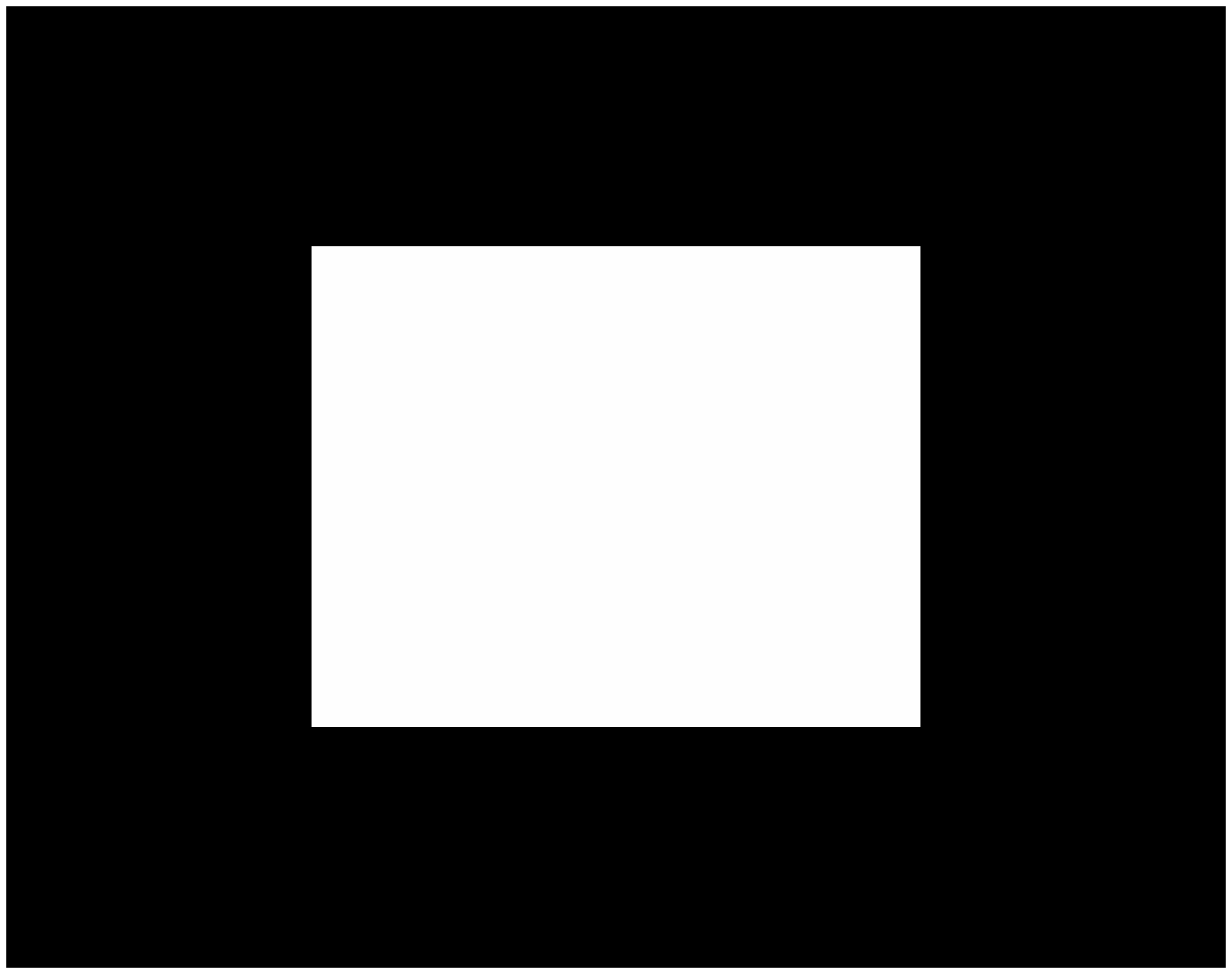}
   \fi
\end{center}
\end{figure}

Using the stable approximation $\rasre$ given by \eqr{rasre}, from the initial guess
$\itr{{u}}{0}$ we generate the sequence $\seq{{u}}$ by
\begin{equation}\label{imp}
\itr{{u}}{n+1} = \rasrne \itr{u}{n}
\quad
\mbox{where}
\quad
\betnpo = \bet_0 + (1-\bet_0)\left(1 - \exp\left(-(n/7)^3\right)\right).
\end{equation}
The rule for updating $\betn$ is a smooth approximation to a step function
from the value $\bet_0$ to the value $1$ centered at iteration $n=7$. 
We compare this algorithm to the HIO \eqr{e:hio} and HPR \eqr{e:hpr} 
algorithms using the same stable projection approximation.  
We study algorithm performance with noisy
data.  The initial points ${u}_0$ are chosen to be the normalized
characteristic function of the support constraint shown in \figr{true}(c).

The data consists of the support/nonnegativity constraint, shown in
\figr{true}(c), and Fourier magnitude data $m$, shown in \figr{true}(b),
with additive noise $\eta$ -- a symmetric, randomly generated 
array with a zero mean Gaussian distribution. 
The signal-to-noise ratio (SNR) is
$20\log_{10}\|m\|/\|u\|=34$~dB.  As motivated in \cite{BCL2}, 
the error metric we use to monitor the algorithms, $\esp$,  is given by
\begin{equation}\label{e:bmeasure}
\espn =
\frac{\big\| \PSP\big(\PM (\itr{{u}}{n})\big)-\PM (\itr{{u}}{n})\big\|^2}
{\big\|\PM (\itr{{u}}{n})\big\|^2}.
\end{equation}
We compute the mean value of the error measure $E_{S+}$ over 100 trials 
with different realizations of the noise and the same initial guess.

\begin{figure}
\caption{\label{good news}
  {\small Error metric $\espn$ averaged over $100$ realizations of noise
  (SNR=34~dB). For (a)-(c) the relaxation parameter for the respective
  algorithms, $\betn$, is fixed.  For (d) $\betn$ varies from $0.75$ to $1.0$
  according to \eqr{imp}.
}}
\begin{center}
  \ifpdf
	  (a) \includegraphics[height=6.0cm,width=6.0cm,angle=0]%
	  {pictures/rasr/l64c_n02_b75_dBmean_100it.png}
	  (b) \includegraphics[height=6.0cm,width=6.0cm,angle=0]%
	  {pictures/rasr/l64c_n02_b87_dBmean_100it.png}\\
	  (c) \includegraphics[height=6.0cm,width=6.0cm,angle=0]%
	  {pictures/rasr/l64c_n02_b99_dBmean_100it.png}
	  (d) \includegraphics[height=6.0cm,width=6.0cm,angle=0]%
	  {pictures/rasr/l64c_n02_vb75_dBmean_100it.png}
	\else
  		(a) \includegraphics[height=6.0cm,width=6.0cm,angle=0]%
		{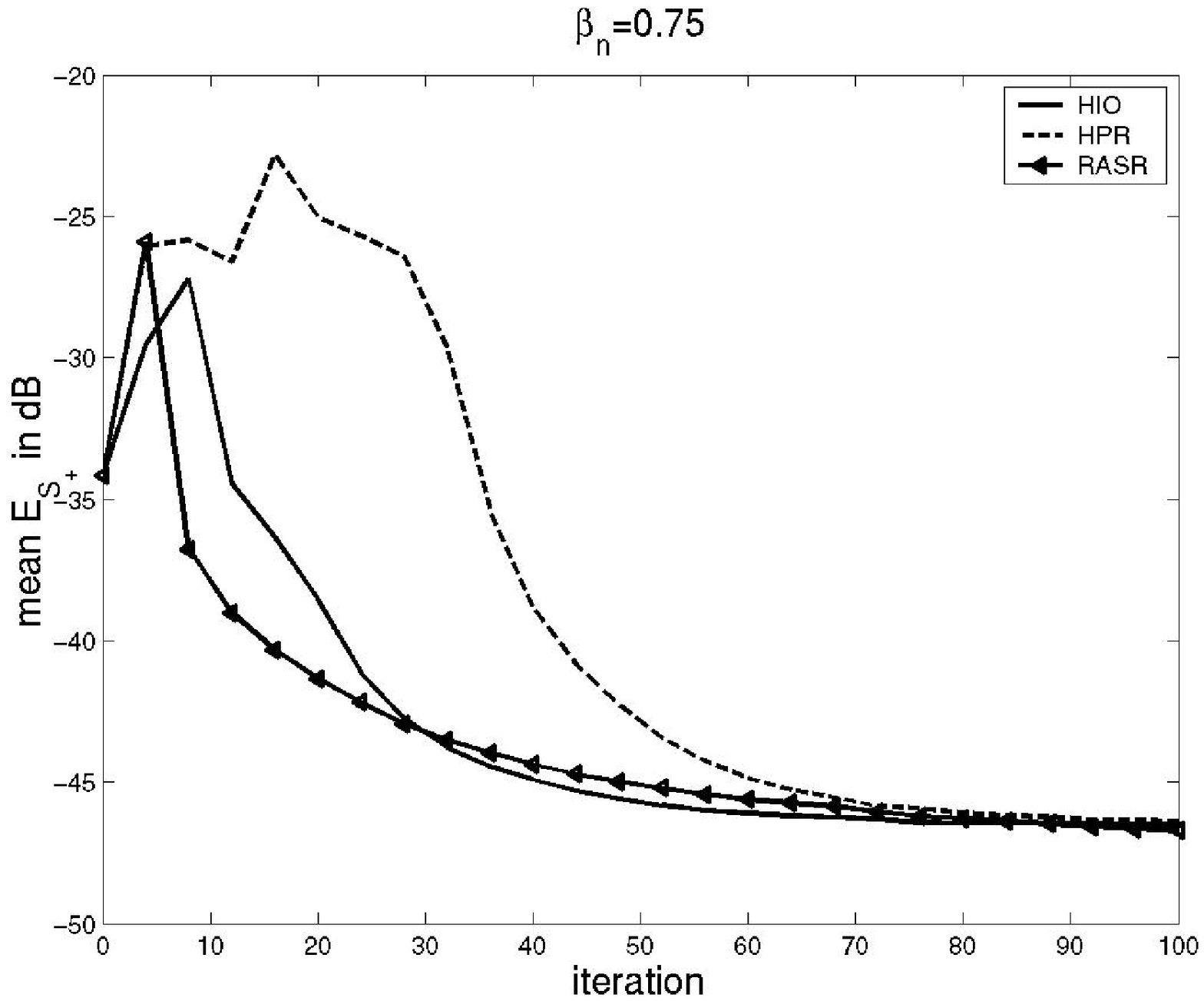}
  		(b) \includegraphics[height=6.0cm,width=6.0cm,angle=0]%
      {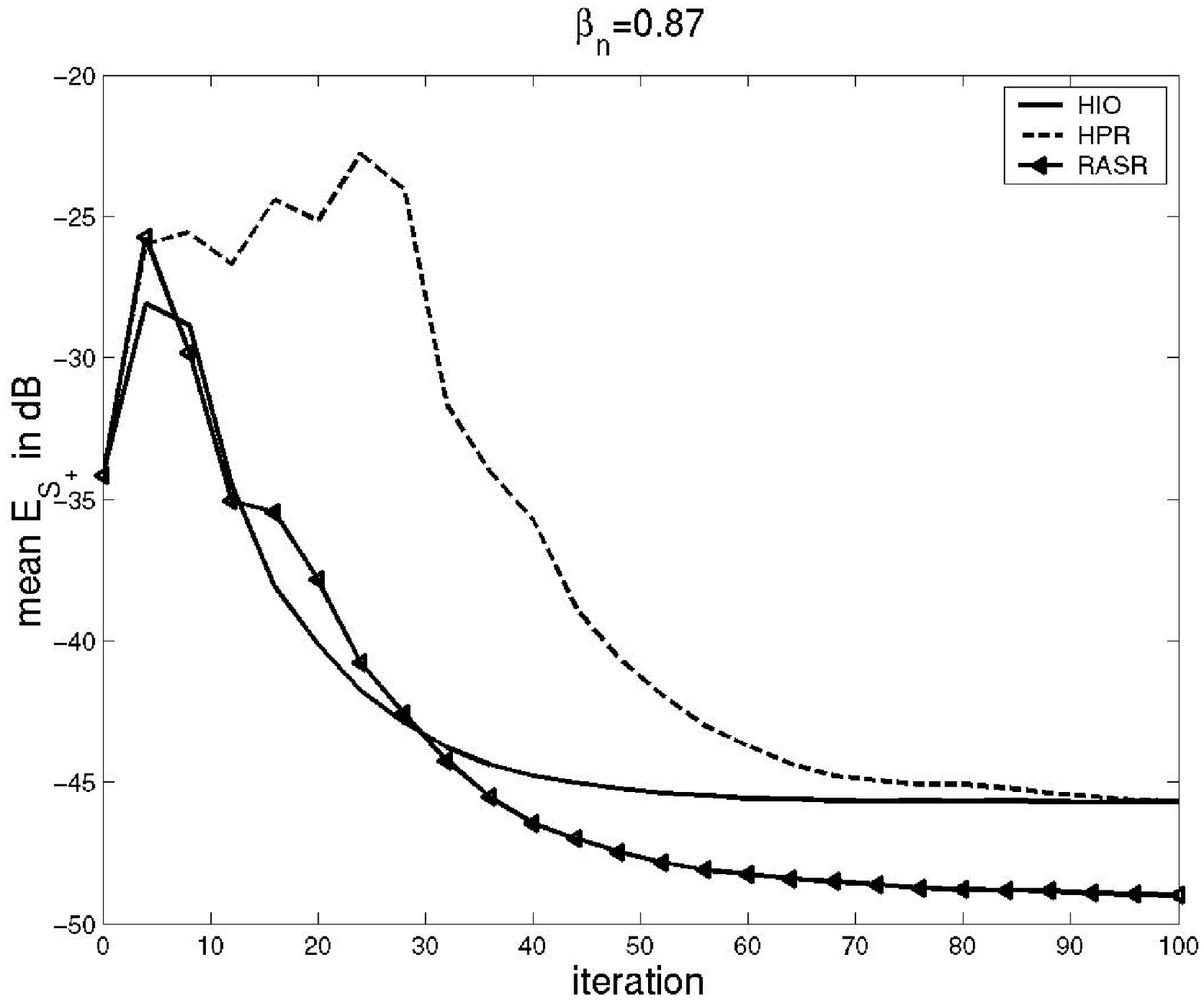}\\
 	   (c) \includegraphics[height=6.0cm,width=6.0cm,angle=0]%
 	   {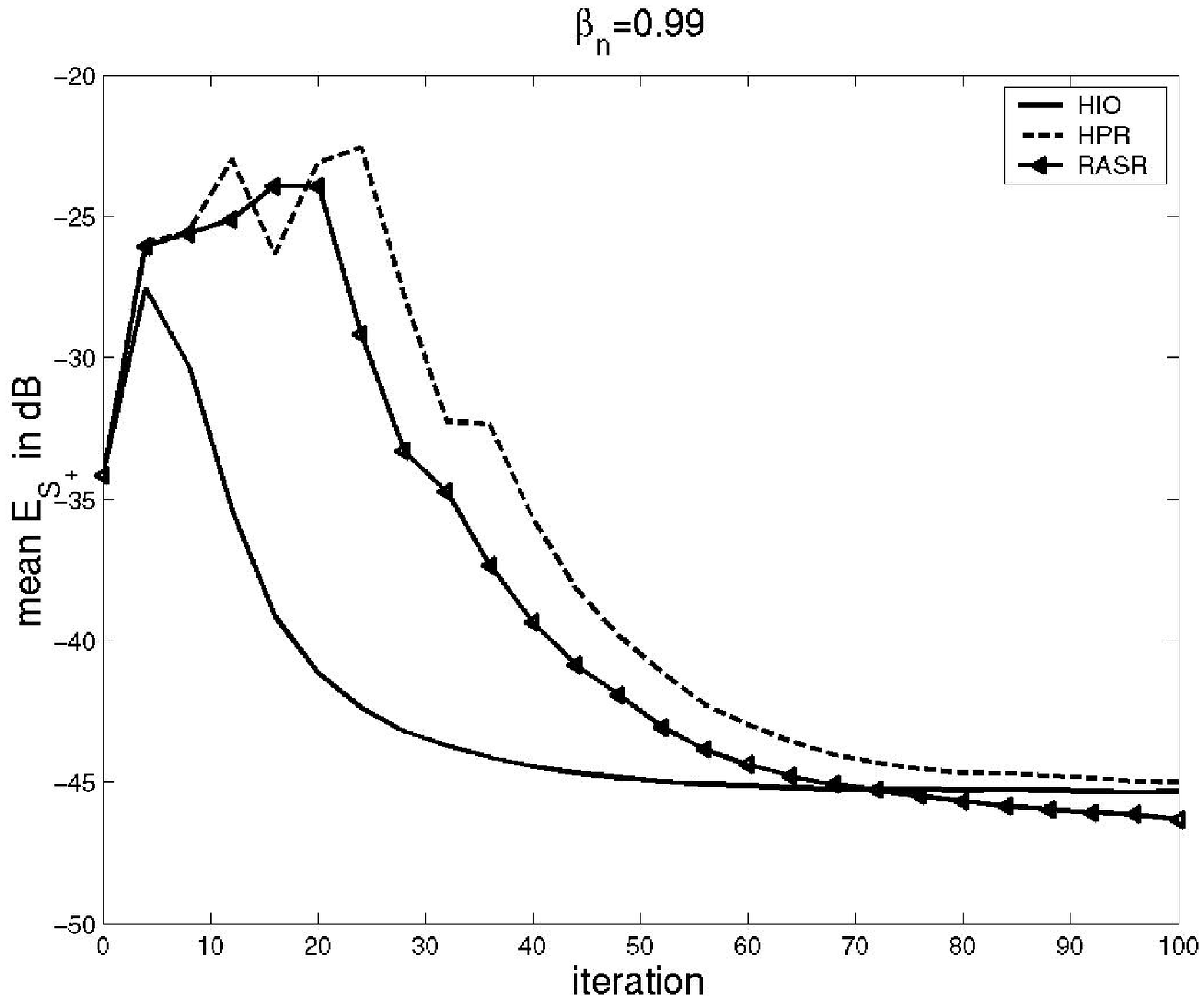}
	   (d) \includegraphics[height=6.0cm,width=6.0cm,angle=0]%
	   {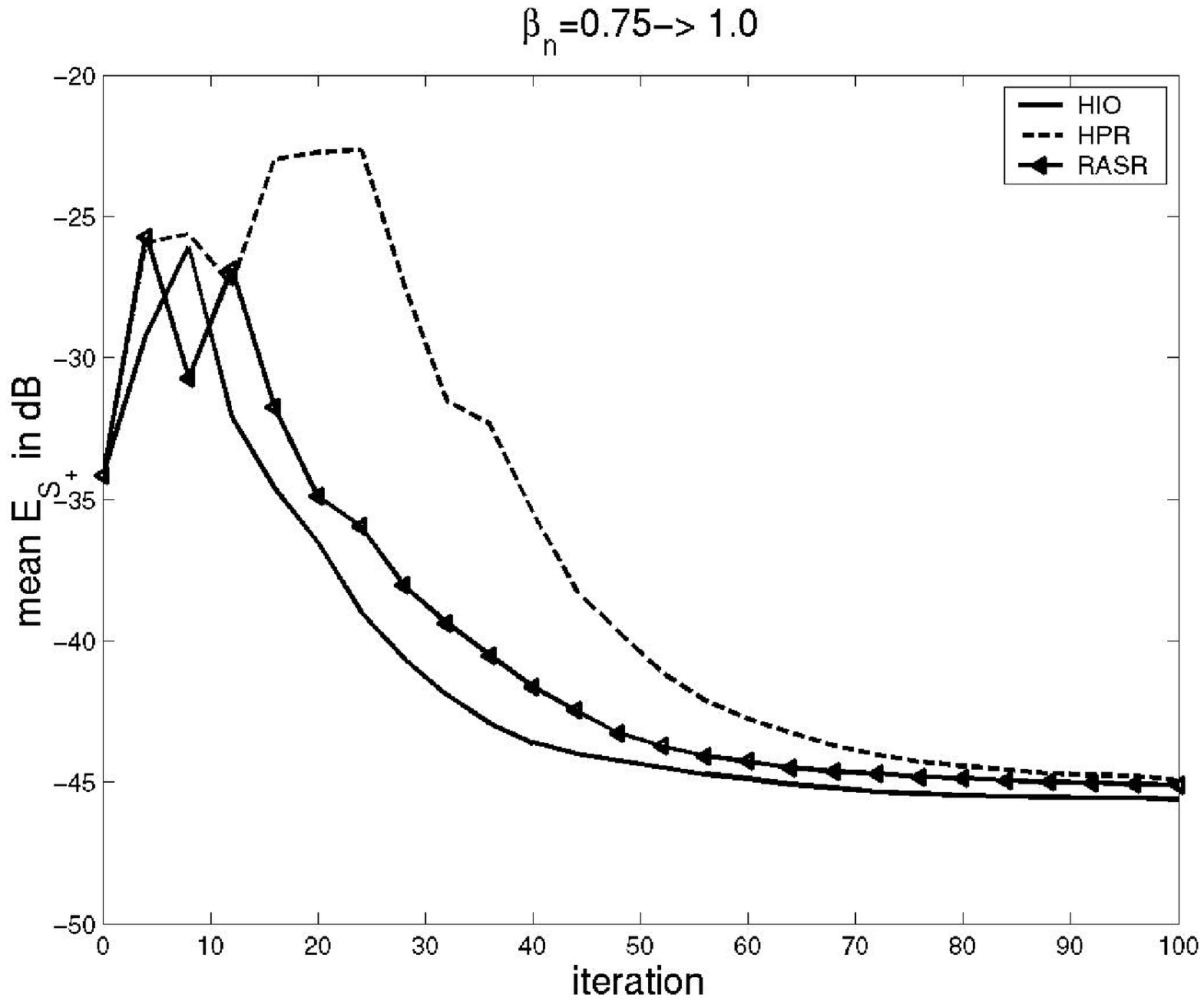}
	\fi
\end{center}
\end{figure}

\begin{figure}
\caption{\label{bad news}
{\small
Typical images recovered after $35$ iterations of the HIO,
HPR, and RAAR algorithms for different relaxation strategies with the
same realization of data noise
(SNR=34~dB) and the same normalized initial guess.  The
variable $\betn$ trials were generated according to the rule given by \eqr{imp}.
}}
\begin{center}
\hskip3.1cm  HIO\hskip2.2cm HPR\hskip2.1cm RAAR
\end{center}
\begin{center}
	\ifpdf
		$\beta=0.75$\hskip1.45cm
		\includegraphics[height=3.0cm,width=3.0cm,angle=0]%
		{pictures/rasr/HIO75_l64c_n02_35it_seed25.jpg}
		\includegraphics[height=3.0cm,width=3.0cm,angle=0]%
		{pictures/rasr/HPR75_l64c_n02_35it_seed25.jpg}
		\includegraphics[height=3.0cm,width=3.0cm,angle=0]%
		{pictures/rasr/RASR75_l64c_n02_35it_seed25.jpg}
	\else
		$\beta=0.75$\hskip1.45cm
		\includegraphics[height=3.0cm,width=3.0cm,angle=0]%
		{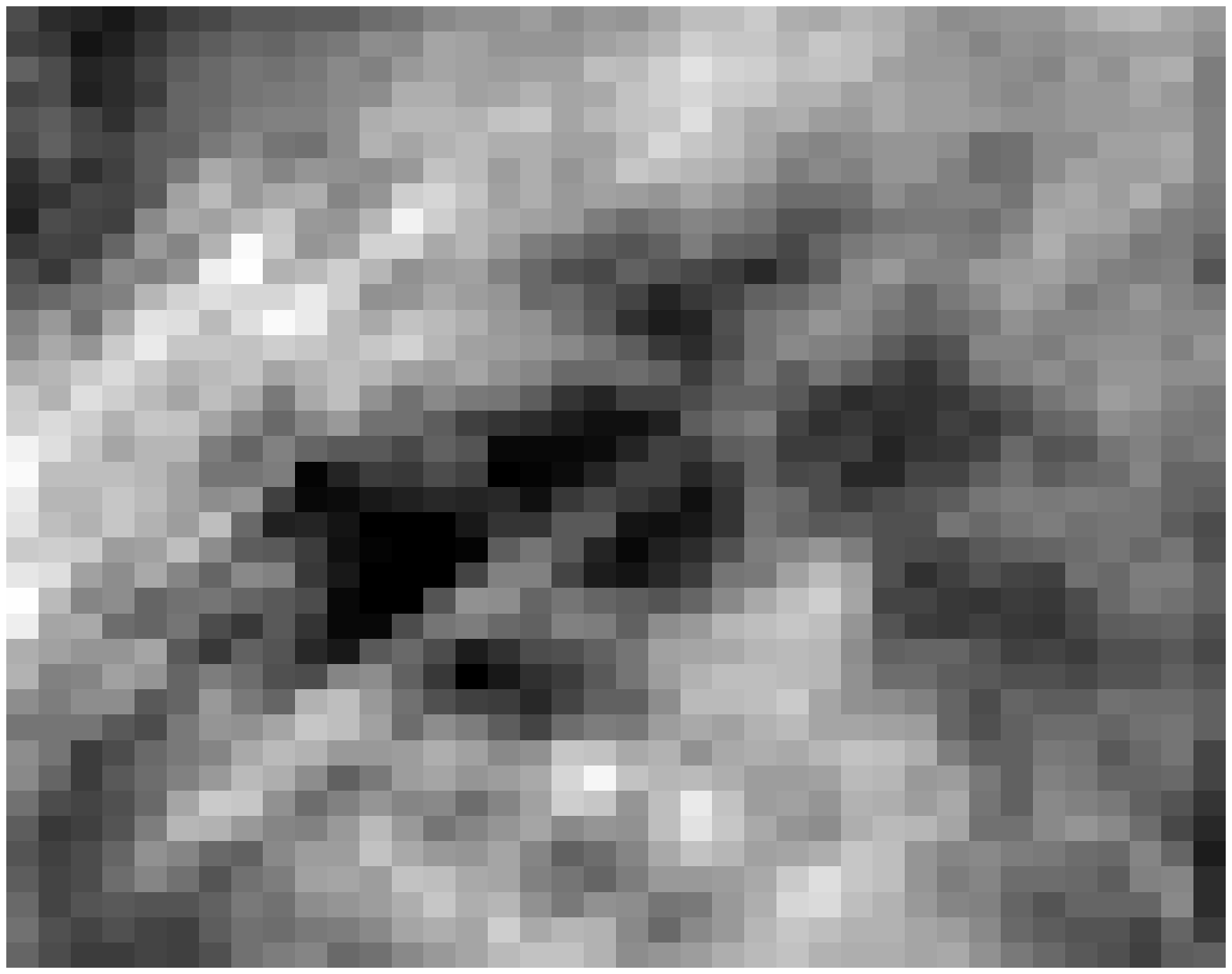}
		\includegraphics[height=3.0cm,width=3.0cm,angle=0]%
		{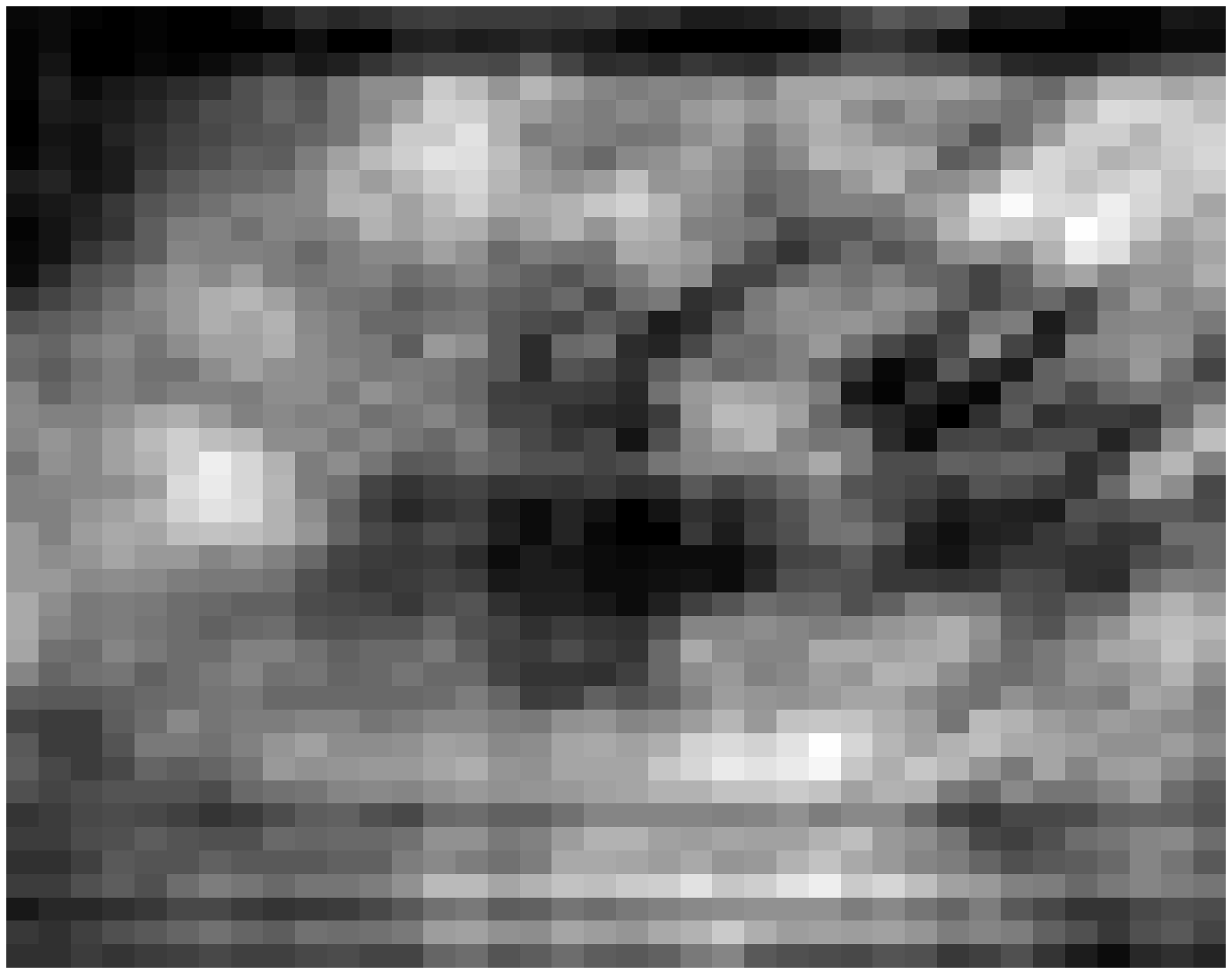}
		\includegraphics[height=3.0cm,width=3.0cm,angle=0]%
		{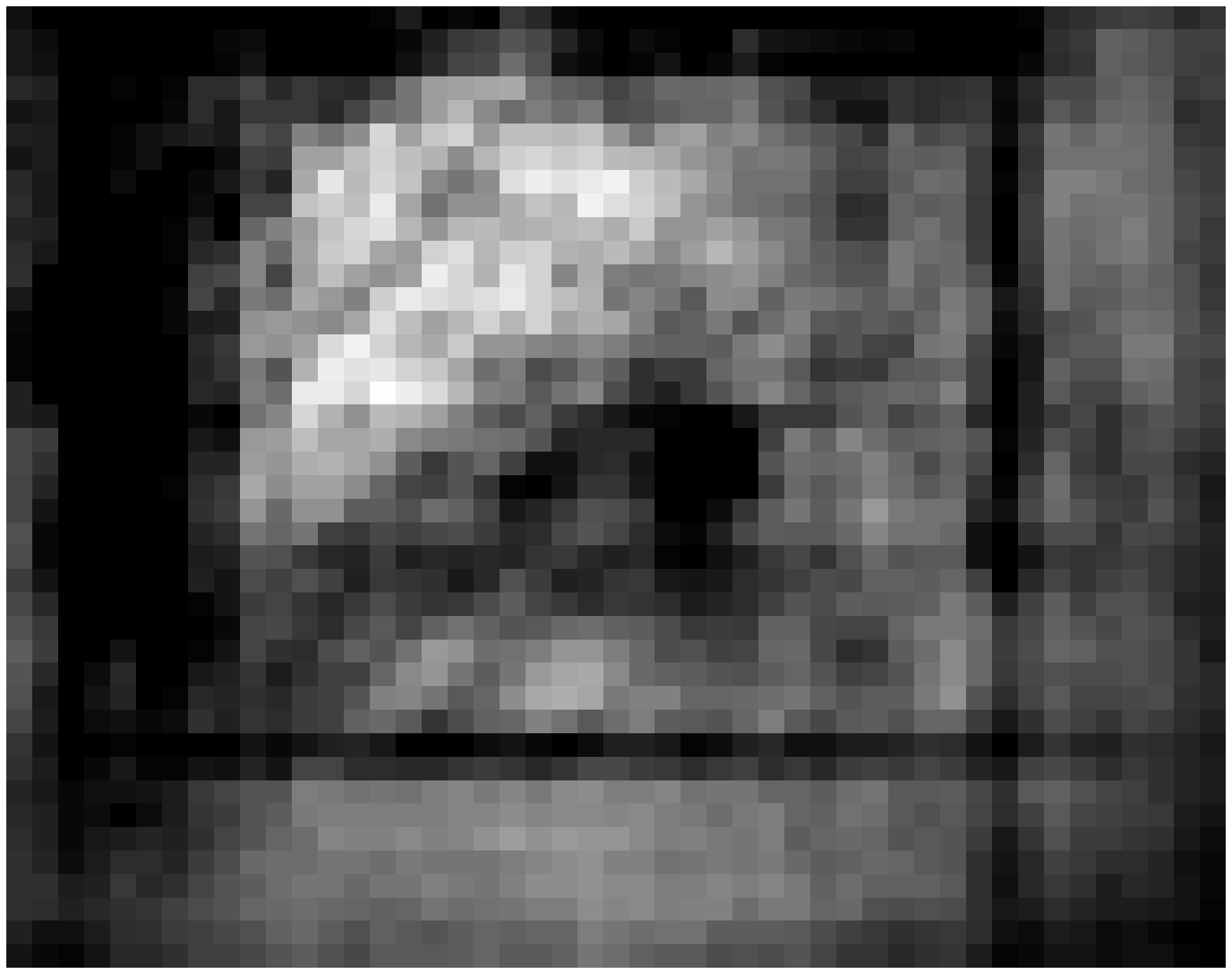}		
	\fi
\end{center}
% \begin{center}
% $\beta=0.87$\hskip1.45cm
% % \epsfig{file=pictures/rasr/HIO87_l64c_n02_35it_seed25.eps,height=3.0cm, width=3.0cm}
% \includegraphics[height=3.0cm,width=3.0cm,angle=0]%
% {pictures/rasr/HIO87_l64c_n02_35it_seed25.jpg}
% % {pictures/rasr/HIO87_l64c_n02_35it_seed25.eps}
% % \epsfig{file=pictures/rasr/HPR87_l64c_n02_35it_seed25.eps,height=3.0cm, width=3.0cm}
% \includegraphics[height=3.0cm,width=3.0cm,angle=0]%
% {pictures/rasr/HPR87_l64c_n02_35it_seed25.jpg}
% % {pictures/rasr/HPR87_l64c_n02_35it_seed25.eps}
% % \epsfig{file=pictures/rasr/RASR87_l64c_n02_35it_seed25.eps,height=3.0cm, width=3.0cm}
% \includegraphics[height=3.0cm,width=3.0cm,angle=0]%
% {pictures/rasr/RASR87_l64c_n02_35it_seed25.jpg}
% % {pictures/rasr/RASR87_l64c_n02_35it_seed25.eps}
% \end{center}
\begin{center}
	\ifpdf
		$\beta=0.87$\hskip1.45cm
		\includegraphics[height=3.0cm,width=3.0cm,angle=0]%
		{pictures/rasr/HIO87_l64c_n02_35it_seed25.jpg}
		\includegraphics[height=3.0cm,width=3.0cm,angle=0]%
		{pictures/rasr/HPR87_l64c_n02_35it_seed25.jpg}
		\includegraphics[height=3.0cm,width=3.0cm,angle=0]%
		{pictures/rasr/RASR87_l64c_n02_35it_seed25.jpg}
	\else
		$\beta=0.87$\hskip1.45cm
		\includegraphics[height=3.0cm,width=3.0cm,angle=0]%
		{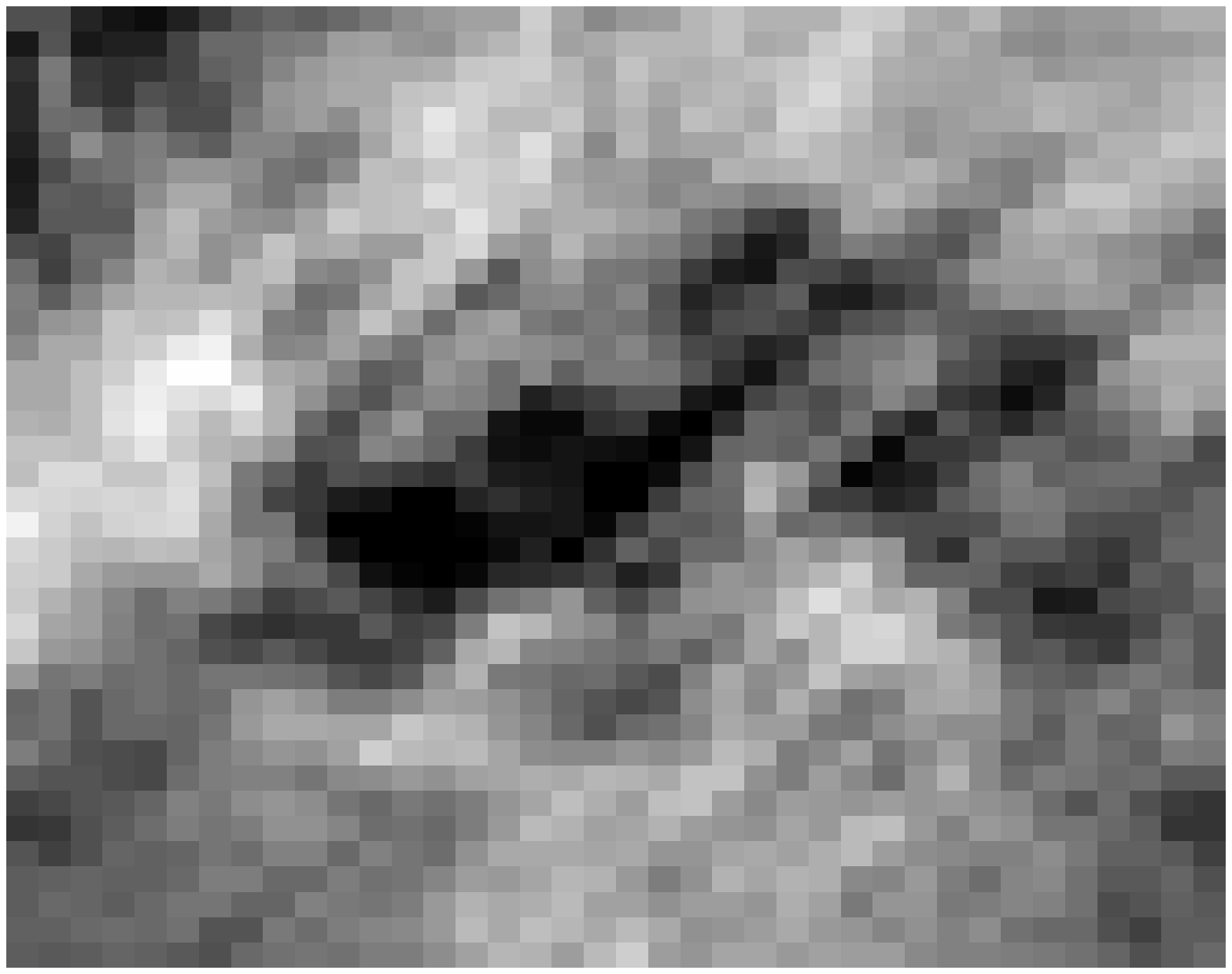}
		\includegraphics[height=3.0cm,width=3.0cm,angle=0]%
		{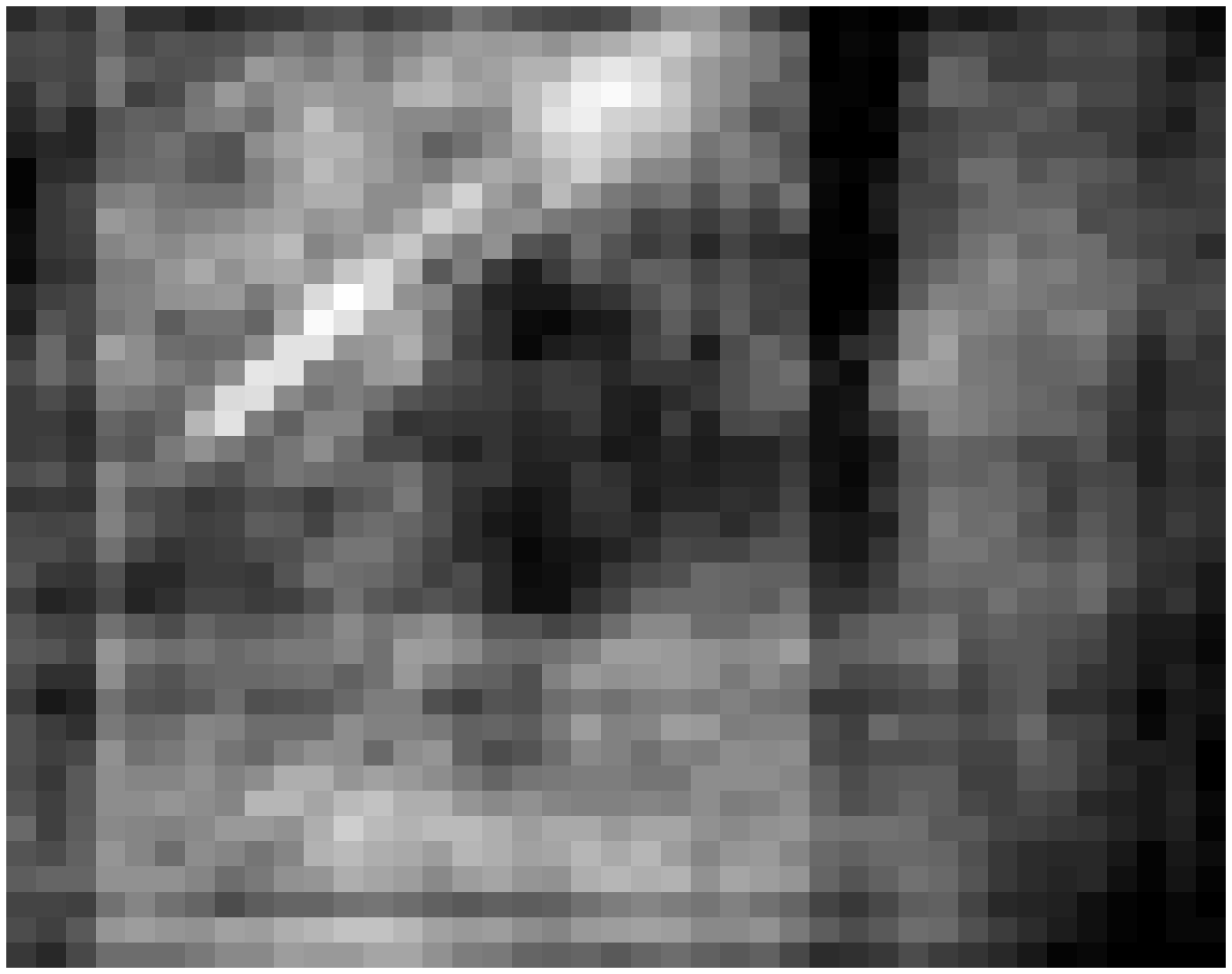}
		\includegraphics[height=3.0cm,width=3.0cm,angle=0]%
		{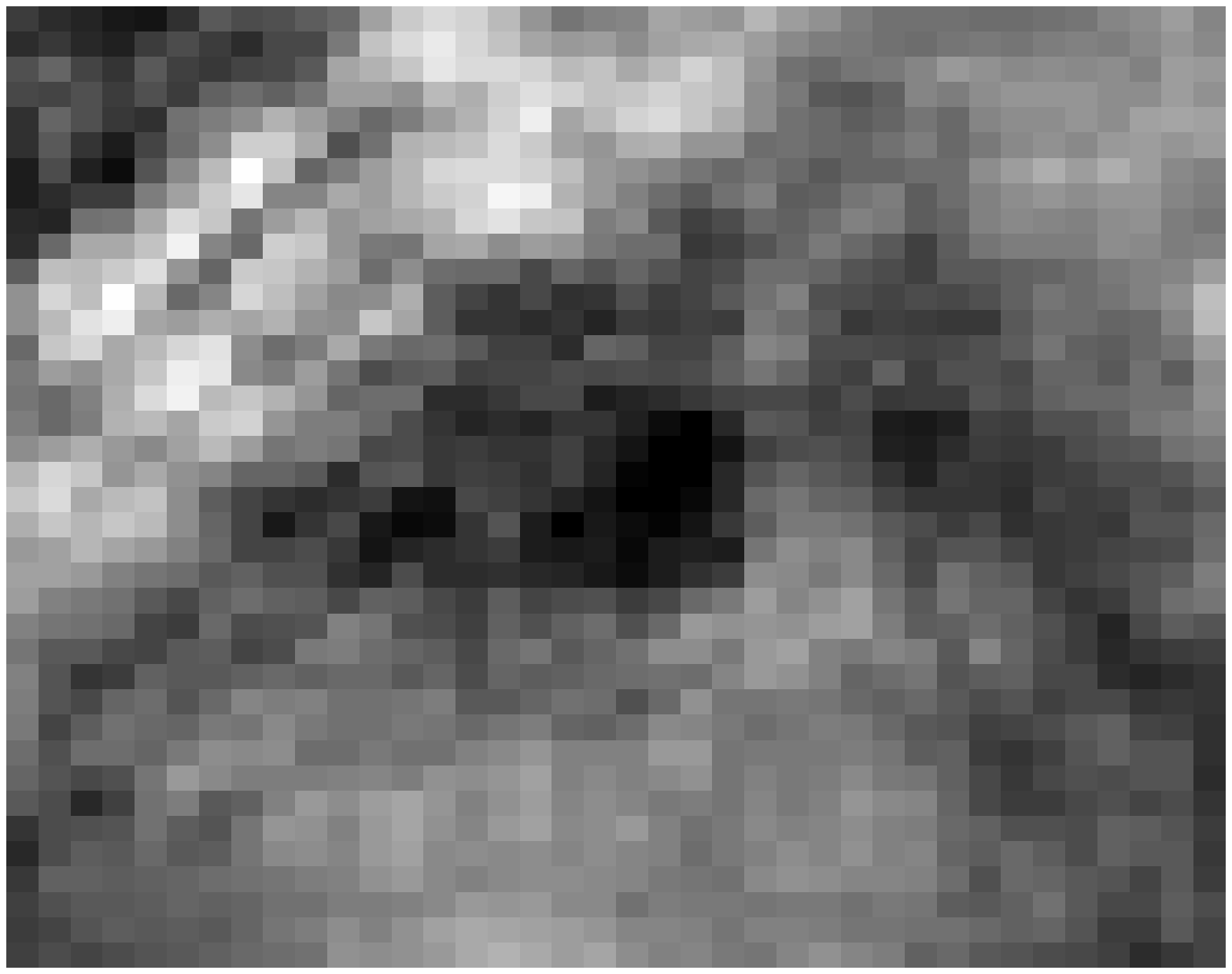}
	\fi
\end{center}
\begin{center}
	\ifpdf
		$\beta=0.99$\hskip1.45cm
		\includegraphics[height=3.0cm,width=3.0cm,angle=0]%
		{pictures/rasr/HIO99_l64c_n02_35it_seed25.jpg}
		\includegraphics[height=3.0cm,width=3.0cm,angle=0]%
		{pictures/rasr/HPR99_l64c_n02_35it_seed25.jpg}
		\includegraphics[height=3.0cm,width=3.0cm,angle=0]%
		{pictures/rasr/RASR99_l64c_n02_35it_seed25.jpg}
	\else
		$\beta=0.99$\hskip1.45cm
		\includegraphics[height=3.0cm,width=3.0cm,angle=0]%
		{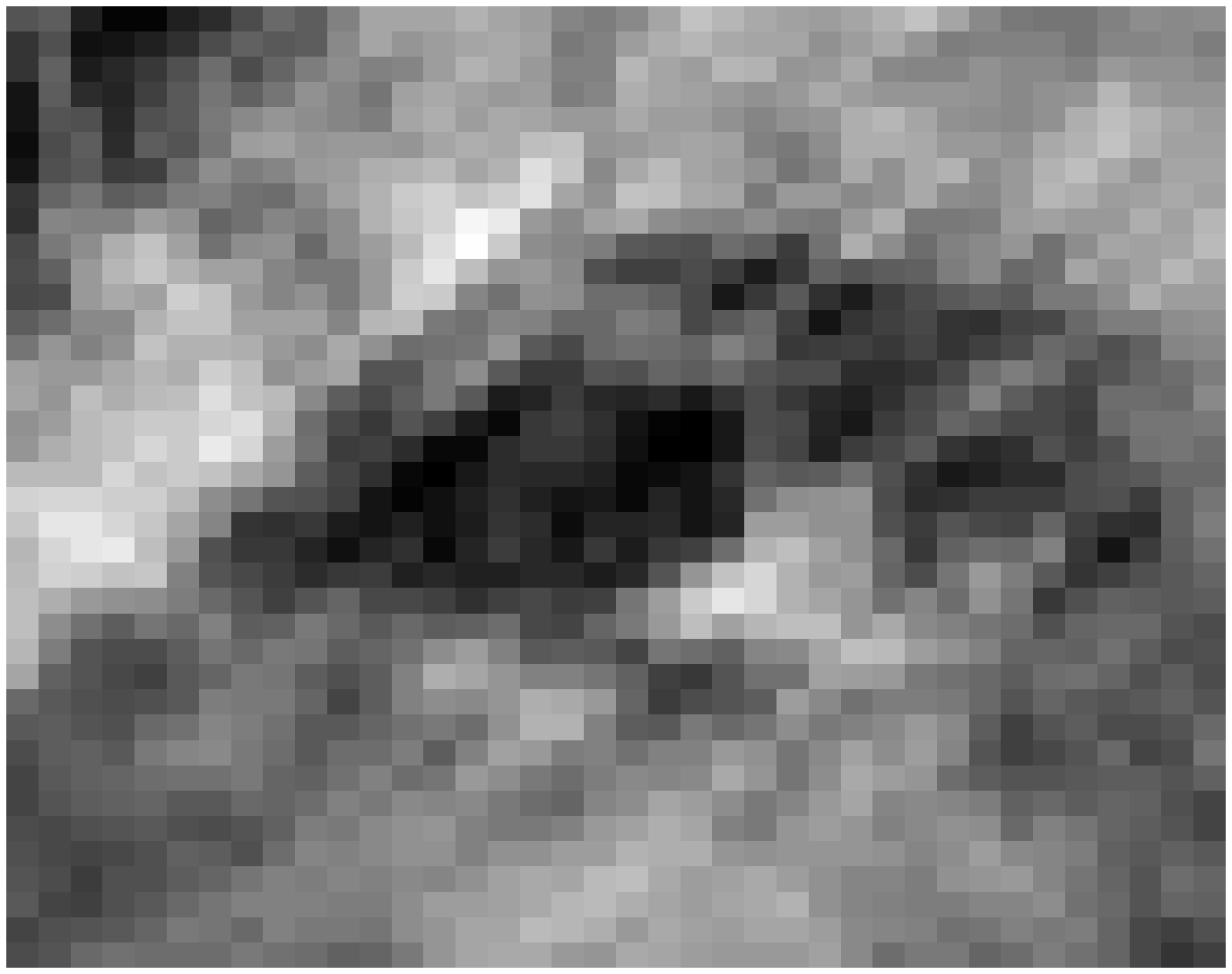}
		\includegraphics[height=3.0cm,width=3.0cm,angle=0]%
		{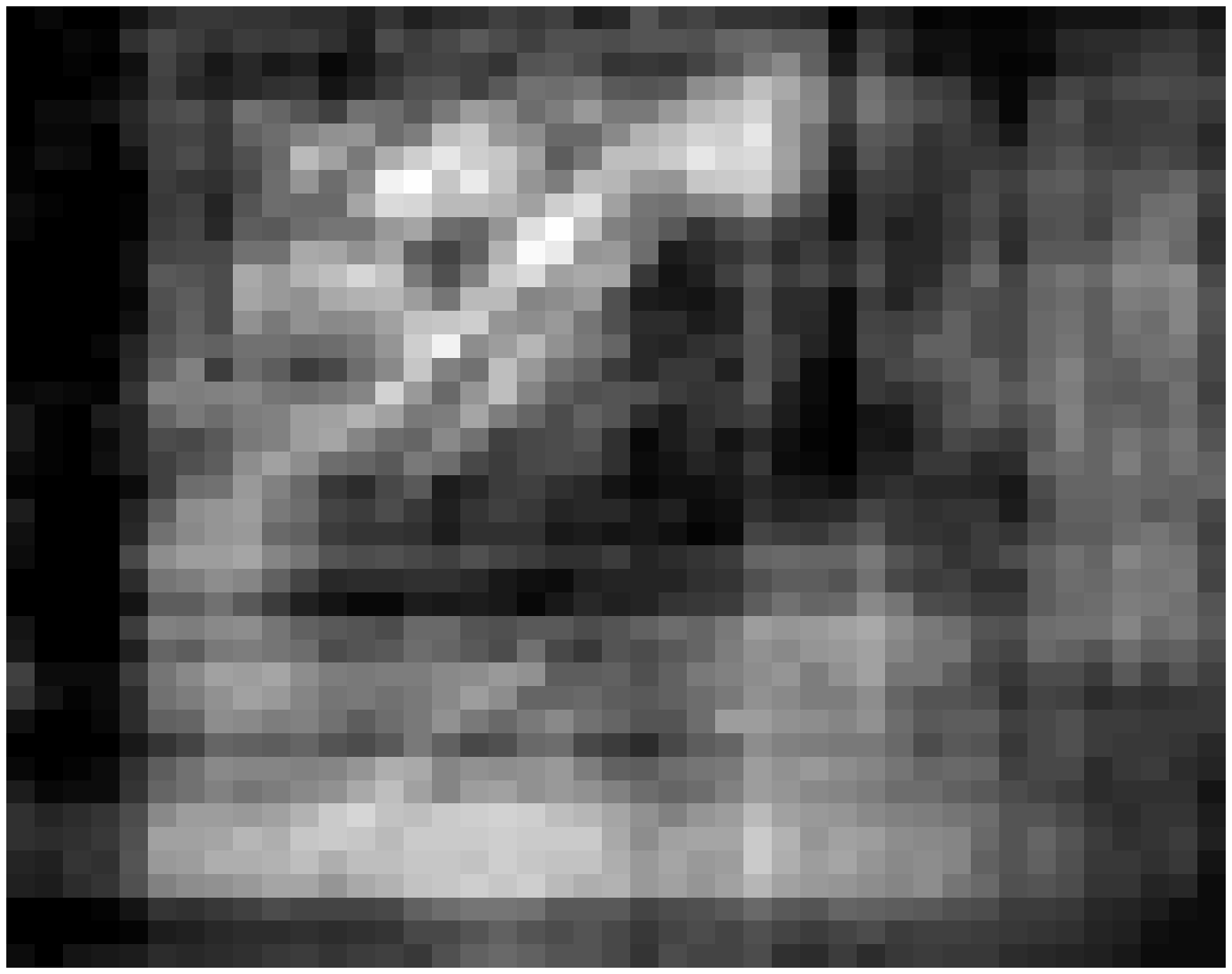}
		\includegraphics[height=3.0cm,width=3.0cm,angle=0]%
		{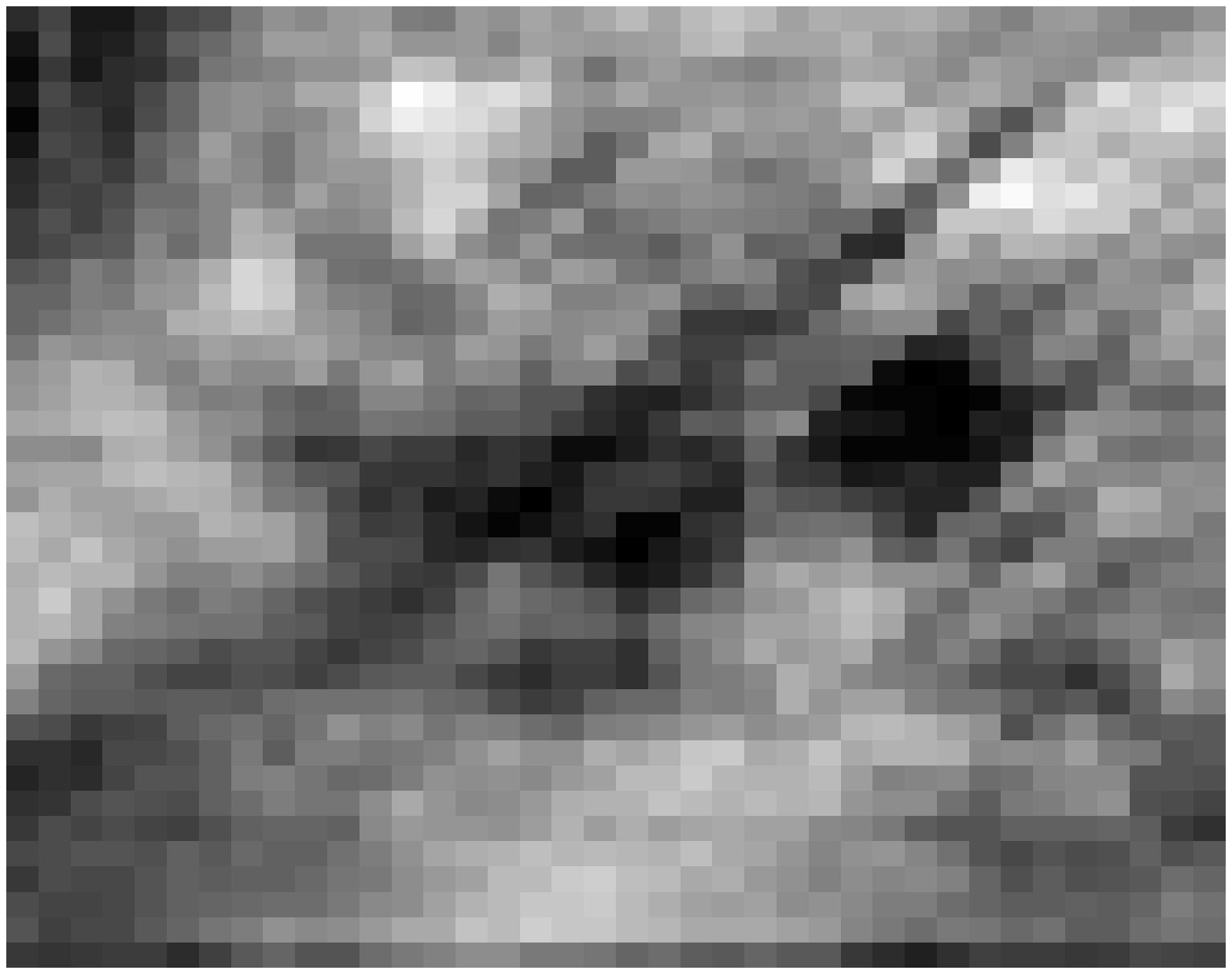}
	\fi
\end{center}
\begin{center}
	\ifpdf
		$\beta=0.75\to 1.0$\hskip.25cm
		\includegraphics[height=3.0cm,width=3.0cm,angle=0]%
		{pictures/rasr/HIOvb75_l64c_n02_35it_seed25.jpg}
		\includegraphics[height=3.0cm,width=3.0cm,angle=0]%
		{pictures/rasr/HPRvb75_l64c_n02_35it_seed25.jpg}
		\includegraphics[height=3.0cm,width=3.0cm,angle=0]%
		{pictures/rasr/RASRvb75_l64c_n02_35it_seed25.jpg}	
	\else
		$\beta=0.75\to 1.0$\hskip.25cm
		\includegraphics[height=3.0cm,width=3.0cm,angle=0]%
		{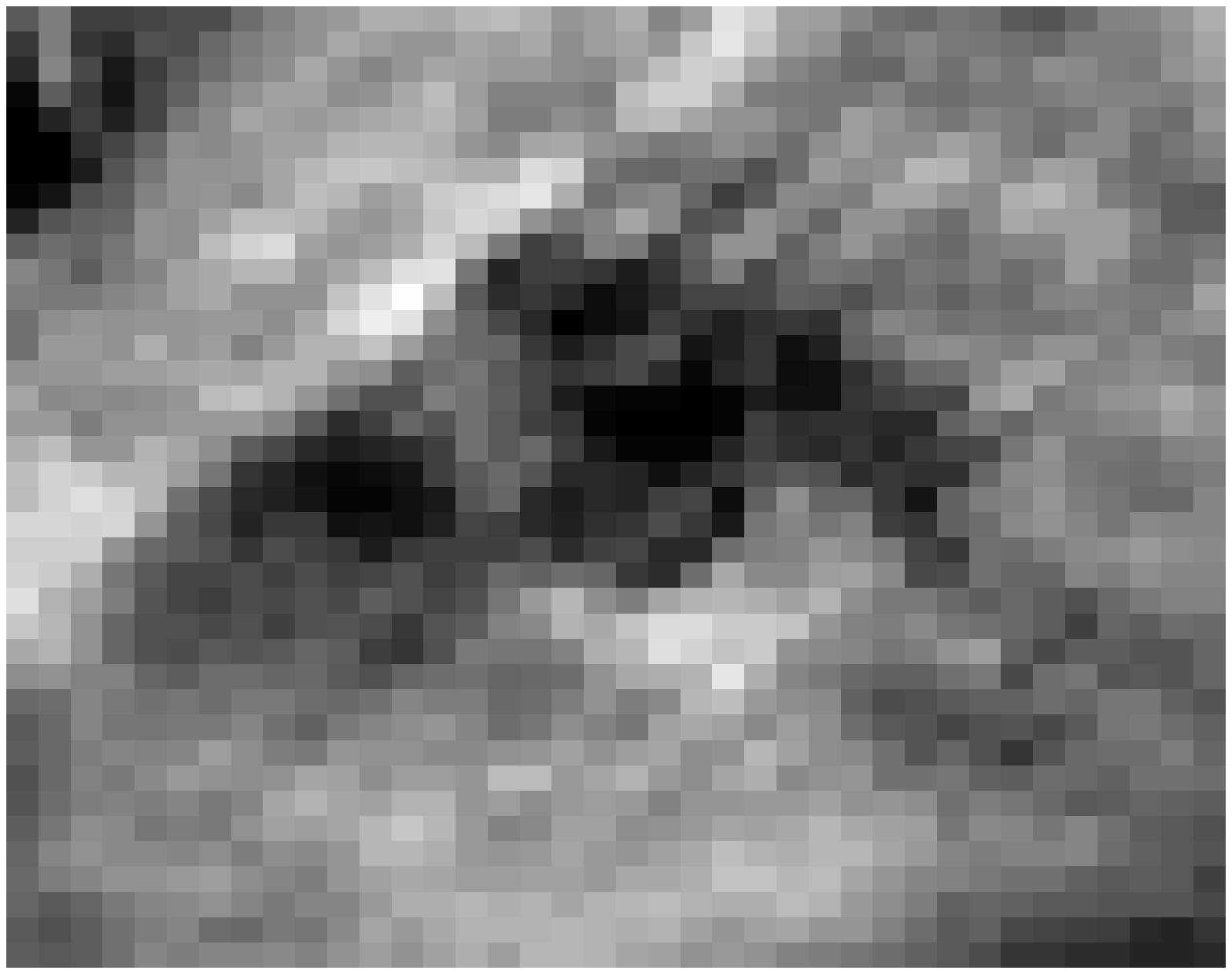}
		\includegraphics[height=3.0cm,width=3.0cm,angle=0]%
		{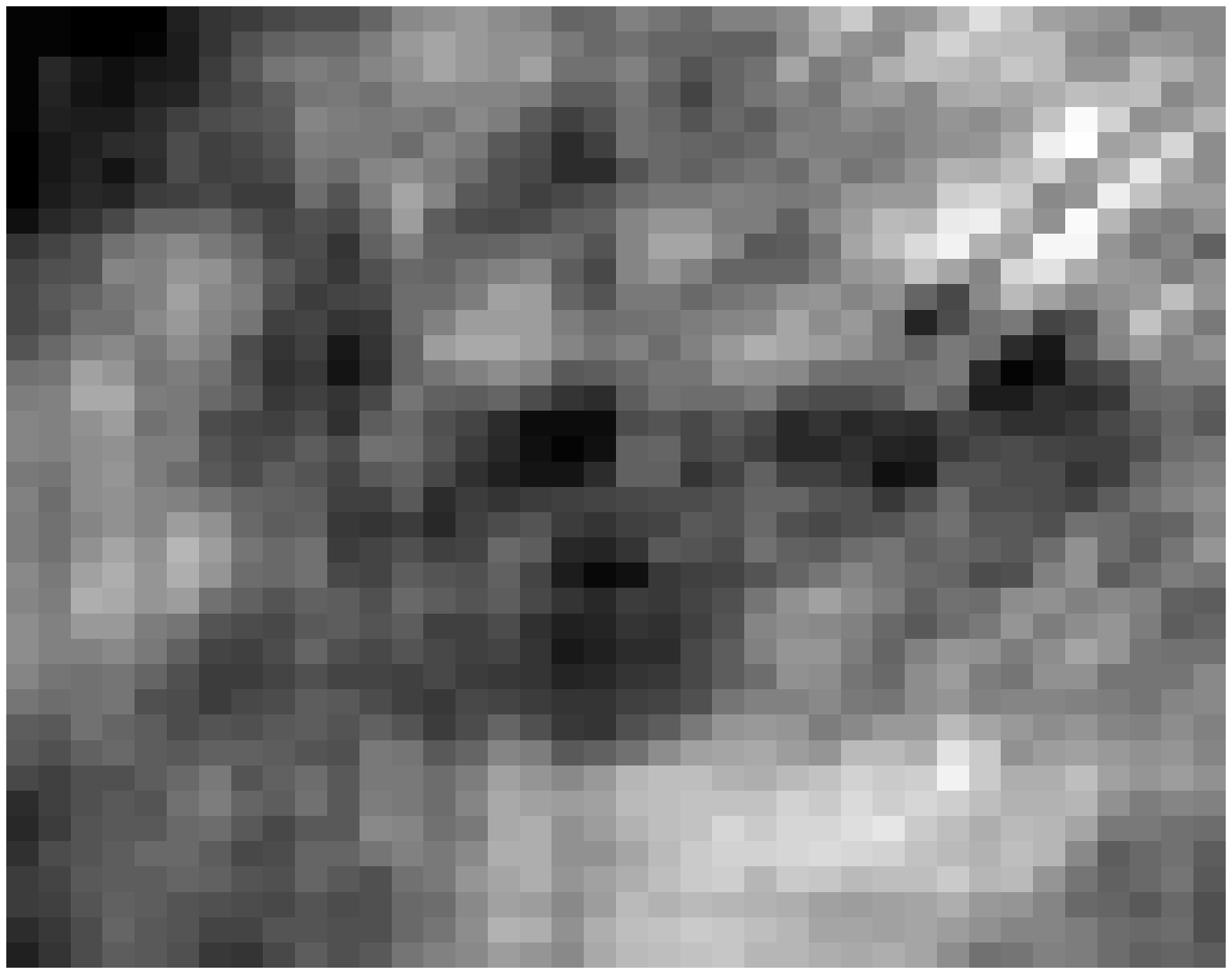}
		\includegraphics[height=3.0cm,width=3.0cm,angle=0]%
		{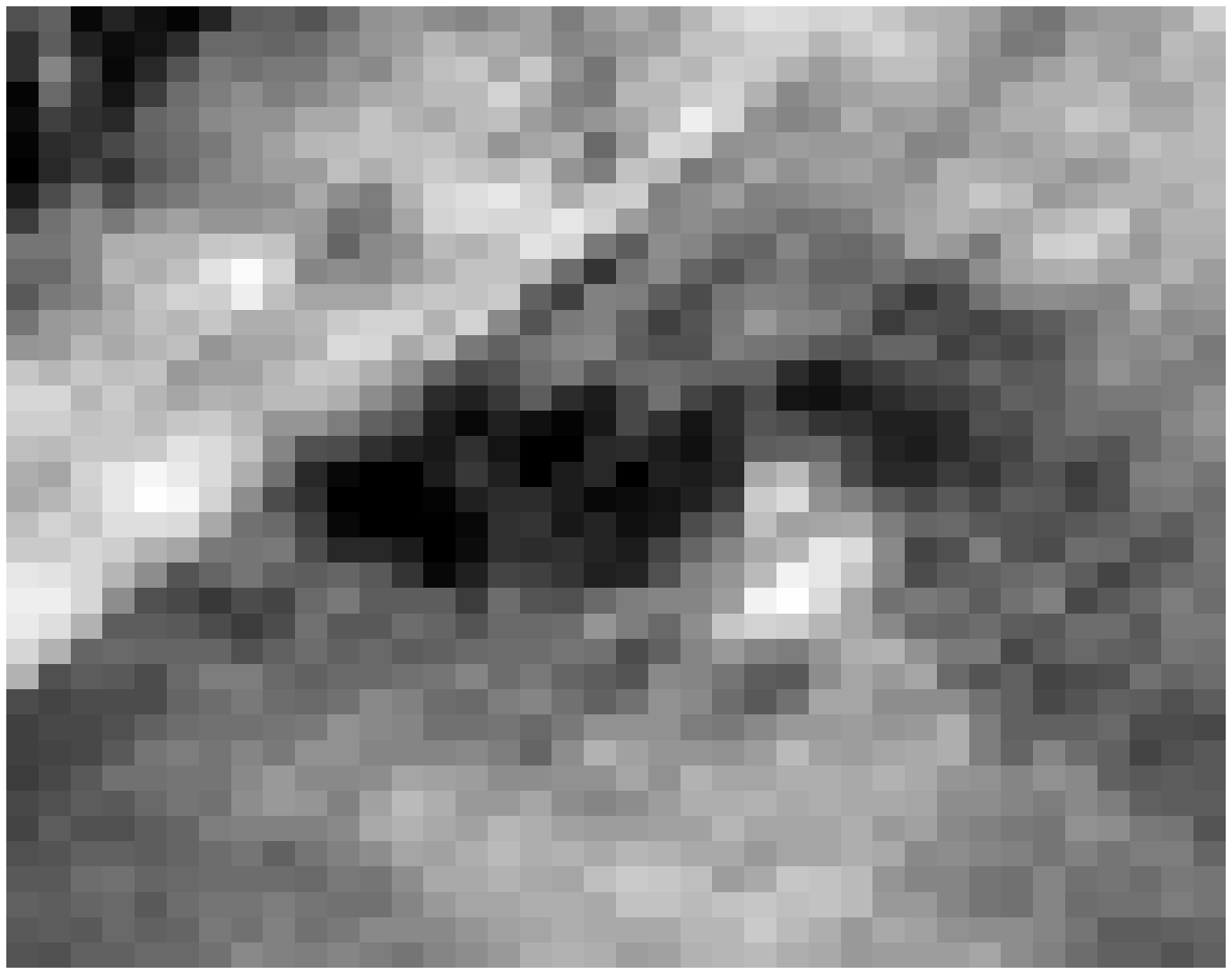}
	\fi
\end{center}
\end{figure}

First, we compare the mean behavior over 100 iterations of two
sets of realizations of the algorithms, each corresponding to
different relaxation strategies, $\beta=0.75,~ \beta=0.87,~ \beta=0.99$
and variable $\betn$ governed by \eqr{imp} with $\bet_0=0.75$.   The
average value of the error
metric at iteration $n$, $\espn$, is shown in Figure~\ref{good news}.
These are all given in decibels (recall that the decibel value of $\alpha>0$ is
$10\log_{10}(\alpha)$).  In \figr{bad news} we show typical estimates
generated by the respective algorithms at iteration $35$,
all from the same realization of noise and the same initial guess.
While the RAAR algorithm with $\beta=0.75$ appears to
perform well as measured by $\esp$ (see \figr{good news}(a)), it is
clear from \figr{bad news} that the quality of solutions found by the RAAR
algorithm degrades rapidly as the relaxation parameter $\beta$ becomes 
small.  For values of $\beta$ near $1.0$ the quality of the iterates 
generated by the RAAR algorithm does 
eventually improve, however, as with the HPR algorithm, it takes many
more iterations to achieve this imporvement.  For static values of $\beta$
the best performance for the RAAR algorithm appears to be achieved with a 
value of $\beta=0.87$.  The variable $\betn$ trials for the RAAR
algorithm yielded the best overall results, measured both by the error
metric, as well as observed picture quality.  
In contrast to this, the relaxation parameter does 
not appear to have any identifiable effect on the performance of the HIO
or HPR algorithms.  

\section{Concluding Remarks}

There are infinitely many relaxation strategies one could implement 
for iterative transform methods, but very few of them admit a meaningful
mathematical analysis.  The standard for phase retrieval algorithms, Fienup's
HIO algorithm, has been identified in a special case with the promising HPR algorithm, 
which in turn, has been identified as a special case of Elser's difference map.
For each of these algorithmic frameworks, the mathematical properties of the
algorithms vary drammatically with the parameter values in a manner analogous to 
bifucations of dynamical systems.  A complete mathematical analysis must
treat all relevant intervals of parameter values on a case by case basis.  No 
such analysis is available for the HIO, HPR or difference map algorithms.  To 
circumvent these difficulties and to improve upon the HPR 
algorithm, we propose a simple relaxation, the RAAR algorithm, of a well 
understood Averaged Averaged Reflection (AAR) algorithm.  The relaxation
is a convex combination of the AAR fixed point operator, and the 
projection onto the {\em data}. This intuitive framework is mathematically
tractable and provides an easy strategy for the choice of relaxation
parameter that, moreover, improves algorithm performance.  
In contrast, it appears that similar relaxation
strategies have little effect on either the HIO or the HPR algorithm. 
We cannot suggest a rule by which to select a 
static value of $\beta$ -- this depends on the data.  
Nevertheless, based on the results for the variable $\betn$ trials, 
we can recommend the fairly generic dynamic relaxation strategy of 
\eqr{imp} for getting the best performance from the RAAR algorithm.
Here the algorithm is significantly relaxed in the early iterations, 
helping the algorithm quickly to find a neighborhood of the solution 
while maintaining fidelity to the data, and then decreasing the
relaxation (i.e. increasing $\betn$) in the neighborhood of the solution to avoid
stagnation at a poor local minimum. To stabilize iterates in the domain of attraction
of a solution, a final fixed value of $\beta$ close to, but less than, $1$, say
$\beta=.99999$ should be chosen.  In a technical point, we also proposed a 
smooth perturbation of the magnitude projector \eqr{gradE} 
to improve the numerical stability of computing the projection onto magnitude 
constraints.

\section*{}
% \ack{Acknowledgements. 
This work was supported by a Post-doctoral 
Fellowship from the Pacific Institute for the Mathematical Sciences.
The author would like to thank Veit Elser for pointing out the 
connection between the HPR algorithm and the difference map.
The author would like to give special thanks to Heinz
Bauschke and Patrick Combettes for their careful reading and 
indispensable comments during the preparation of this work.
% }

% \bibliographystyle{siam}\bibliography{acr}
% \referencelist[acr]
\end{document}